\documentclass[11pt]{amsart}
\usepackage[all]{xy}

\usepackage{amssymb,amscd,amsmath,amsthm}

\theoremstyle{plain}
\newtheorem{thm}{Theorem}[section]
\newtheorem{theorem}[thm]{Theorem}

\newtheorem{prop}[thm]{Proposition}

\theoremstyle{definition}
\newtheorem{defn}[thm]{Definition}
\newtheorem{conj}[thm]{Conjecture}

\newtheorem{prob}[thm]{Problem}
\newtheorem{quest}[thm]{Question}

\theoremstyle{remark}

\newtheorem*{remark*}{Remark}

\newcommand{\cA}{{\mathcal{A}}}
\newcommand{\cB}{{\mathcal{B}}}

\newcommand{\cD}{{\mathcal{D}}}

\newcommand{\cF}{{\mathcal{F}}}
\newcommand{\cG}{{\mathcal{G}}}
\newcommand{\cH}{{\mathcal{H}}}
\newcommand{\cI}{{\mathcal{I}}}

\newcommand{\cK}{{\mathcal{K}}}
\newcommand{\cL}{{\mathcal{L}}}

        \newcommand{\field}[1]{{\mathbb{#1}}}
        \newcommand{\NN}{\field{N}}
        \newcommand{\ZZ}{\field{Z}}
        
        \newcommand{\RR}{\field{R}}
        \newcommand{\CC}{\field{C}}

\newcommand{\supp}{\operatorname{supp}}
\newcommand{\Ker}{\operatorname{Ker}}

\newcommand{\Dom}{\mbox{\rm Dom}}

\newcommand{\Tr}{\operatorname{Tr}}
\newcommand{\tr}{\operatorname{tr}}

\newcommand{\res}{{\rm res}\;}

\begin{document}

\title[Noncommutative spectral geometry of foliations]
{Noncommutative spectral geometry of Riemannian foliations: some
results and open problems}
\author{Yuri A. Kordyukov}
\address{Institute of Mathematics, Russian Academy of Sciences, Ufa,
Russia} \email{yuri@imat.rb.ru}

\thanks{Supported by Russian Foundation of Basic Research
(grant no. 04-01-00190)}

\begin{abstract}
We review some applications of noncommutative geometry to the
study of transverse geometry of Riemannian foliations and discuss
open problems.
\end{abstract}
\maketitle

\section*{Introduction}
The main subject of this paper is the Riemannian geometry of the
leaf space of a compact foliated manifold. Moreover, we will
mostly consider the simplest case of the leaf space of a
Riemannian foliation. Our purpose is to explain some basic ideas
and results in noncommutative geometry and its applications to the
study of the leaf space of a foliation and present some open
problems in analysis and geometry on foliated manifolds motivated
by these investigations.

Applications of noncommutative geometry to the study of singular
geometrical objects such as the leaf space of a foliated manifold
are based on several fundamental ideas.

The first idea is to pass from geometric spaces to (analogues of)
algebras of functions on these spaces and translate basic concepts
and constructions to the algebraic language. This is well-known
and has been used for a long time, for instance, in algebraic
geometry.

The second idea is that, in many important cases, it is natural to
consider analogues of algebras of functions on a singular
geometric space to be noncommutative algebras. In
Section~\ref{s:ncspaces}, we describe the construction of
noncommutative algebras associated with the leaf space of a
foliation due to Connes \cite{Co79}. Actually, an arbitrary
noncommutative algebra can be viewed in many cases as an algebra
of functions on some virtual geometric space or, in other words,
as a noncommutative space. For instance, a $C^*$-algebra is the
algebra of continuous functions on a virtual topological space, a
von Neumann algebra is the algebra of essentially bounded
measurable functions on a virtual measurable space and so on.
Therefore, the theory of $C^*$-algebras is a far-reaching
generalization of the theory of topological spaces and is often
called noncommutative topology. The theory of von Neumann algebras
is a generalization of the classical measure and integration
theory and so on. Such a geometric point of view turns out to be
very useful in operator theory and is also well known.

So the correspondence between classical geometric spaces and
commutative algebras is extended to the correspondence between
singular geometric spaces and noncommutative algebras, and we need
to generalize basic concepts and constructions on geometric spaces
to the noncommutative setting. It should be noted that, as a rule,
such noncommutative generalizations are quite nontrivial and have
richer structure and essentially new features than their
commutative analogues.

The main purpose of noncommutative differential geometry, which
was initiated by Connes \cite{Co:nc} and is actively developing at
present time (cf. the recent surveys
\cite{Connes2000,ConnesLNM1831} and the books
\cite{Co,Gracia:book,Landi:book} in regard to different aspects of
noncommutative geometry), is the extension of analysis, the
analytic objects on geometric spaces, to the noncommutative
setting.

We will discuss only one aspect of this theory ---  namely,
Riemannian geometry of singular spaces. Here there is another idea
suggested by Connes: in order to develop Riemannian geometry, one
can start with abstract functional-analytic analogues of natural
geometric operators on a singular space in question and try to
reconstruct basic geometric information from spectral data of
these operators. This idea goes back to spectral geometry.

Usually, spectral geometry is considered as the investigation of a
famous question by Mark Kac: ``Can one hear the shape of a drum?''
If the answer is negative (and now it is known that this is, in
general, so), then the following question is: ``Which geometrical
properties of a drum can one hear?'' We refer the reader, for
instance, to \cite{Berard86,BGM71,Gordon00,Chavel99,Colin96} for
some survey papers on the spectral theory of the Laplace operator
and spectral geometry.

Let $(M,g)$ be a compact Riemannian manifold of dimension $n$,
$\Delta_g$ the associated Laplace-Beltrami operator,
$0=\lambda_0<\lambda_1\leq \lambda_2\leq\lambda_3\leq \cdots,
\lambda_j\to+\infty,$ the set of the eigenvalues of $\Delta_g$
(counted with multiplicities), $\{\varphi_j\in C^\infty(M) :
j=1,2,\ldots \}$ a corresponding complete orthonormal system of
eigenfunctions in $L^2(M)$ such that
$\Delta_g\varphi_j=\lambda_j\varphi_j.$ Consider the eigenvalue
distribution function
\[
N(\lambda)=\sharp\{j: \lambda_j\leq \lambda\}, \quad \lambda \in
\RR.
\]
Recall the following well-known asymptotic formula for
$N(\lambda)$ called the Weyl asymptotic formula:
\[
N(\lambda)=\frac{|B_n|}{(2\pi )^n}{\rm vol}\, M\cdot
\lambda^{n/2}+O(\lambda^{(n-1)/2}), \quad \lambda\to +\infty,
\]
where $|B_n|$ denotes the volume of a unit $n$-dimensional ball.
This formula shows that one can hear the dimension of $M$ and the
volume of $M$. One can also consider the heat trace asymptotic
expansion:
\[
\tr e^{-t\Delta_g}\sim \frac{1}{t^{n/2}} (a_0+ a_1t^{1/2} +a_2t
+\ldots ), \quad t\to 0,
\]
where $a_j$ are integrals of polynomials of the curvature and its
derivatives, or the residues of the zeta-function $\zeta(z)$,
which is defined by the formula
\[
\zeta(z)=\sum_{j=1}^{+\infty}\lambda_j^{-z}, \quad {\rm
Re}\,z>\frac{n}{2},
\]
and extends to a meromorphic function in the entire complex plane.
These formulas allow one to reconstruct some local
differential-geometric invariants from the spectral data of the
Laplace operator.

Among other types of geometric invariants that can be
reconstructed from the spectral data of the Laplace operator
$\Delta_g$, let us mention first the lengths of closed geodesics.
This can be done by considering the singularities of the trace of
the wave group $e^{i t \sqrt{\Delta_g}}$. The
Duistermaat-Guillemin trace formula provides us with more
invariants of the closed geodesics (for instance, so-called wave
invariants and the Birkhoff normal form of the Poincar\'e map),
which can be reconstructed from the spectrum of $\Delta_g$.

To proceed further, we should extend the operator data we are
starting with. First, one can consider the signature operator
$d+d^*$ on differential forms or the Dirac operator on spinors and
use the Hodge theory and the index theory of elliptic operators.
Second, one can take into considerations the algebra of smooth
functions on $M$ considered as an algebra of bounded operators in
$L^2(M)$. This will lead us to local analogues of the facts
mentioned above, say, to the local Weyl asymptotic formula and so
on. Finally, we will arrive at classical mechanical and quantum
mechanical objects on $M$ and relations between these objects
(problems of quantization and semiclassical limits). Let us recall
some basic information on classical and quantum mechanics.

In classical mechanics, a point particle, moving on a compact
manifold $M$ (called the configuration space), is described by a
point of the phase space, which is the cotangent bundle $T^*M$ of
$M$, and the evolution of the phase space point is governed by
Hamilton's equations of motion. In quantum mechanics, a point
particle on a compact manifold $M$ is described by a function in
$L^2(M)$ called the wave function or wave packet. The evolution of
the quantum particle is determined by the Schr\"odinger equation.

In classical mechanics, observables (that is, quantities that we
can observe, such as position, momentum and energy) are
represented by real-valued functions on the phase space. In
quantum mechanics, they are represented by self-adjoint
(unbounded) operators in $L^2(M)$.

In particular, a Riemannian metric $g$ considered as a function on
$T^*M$ is the Hamiltonian (the energy) of a free classical
particle on the configuration space $M$, and the associated
Laplace operator $\Delta_g$ is a Hamiltonian of the free quantum
particle on the configuration space $M$. Therefore, many spectral
quantities we will consider can be treated as quantum analogues
(quantization) of different classical objects, and many classical
objects can be treated as some classical limits. For instance,
quantization of the algebra $C^\infty(M)$ is the subalgebra in
$\cL(L^2(M))$ that consists of the corresponding multiplication
operators. Quantization of the cotangent bundle $T^*M$ is the
algebra of pseudodifferential operators on $M$.

We now extend these ideas to noncommutative algebras. We start
with an involutive algebra $\cA$, a noncommutative analogue of an
algebra of (complex-valued) functions on a singular geometric
object $X$. First, we quantize the algebra $\cA$, taking a
$\ast$-representation of $\cA$ in a Hilbert space $\cH$. Then we
need an abstract analogue $D$ of a first order elliptic
pseudodifferential operator on a compact manifold whose definition
goes back to Atiyah and Kasparov. The resulting object $(\cA, \cH,
D)$ is called a spectral triple or an unbounded Fredholm module
over $\cA$. It can be considered as a virtual (or noncommutative)
geometric space, where $D$ plays the role of a Riemannian metric.
Starting from a spectral triple and using ideas from spectral
geometry, index theory and quantization mentioned above, one can
define analogues of basic geometric and analytic objects on the
associated noncommutative geometric space such as dimension,
differential, differential forms, Riemannian volume form,
cotangent bundle, geodesic flow and so on. A spectral triple can
be associated to a compact Riemannian manifold. In this classical
case, such noncommutative generalizations are shown to be
equivalent to their classical counterparts.

In the case of the leaf space $M/\cF$ of a foliated manifold $(M,
\cF)$, many geometric and analytic objects on this singular space
can be introduced "naively", at the level of sets and points, as
the corresponding holonomy invariant objects on the ambient
manifold. For instance, a holonomy invariant Riemannian metric on
the fibers of the normal bundle of $\cF$ can be considered as a
substitute of a Riemannian metric on $M/\cF$. Such a metric exists
only if the foliation is Riemannian. One can associate a spectral
triple to any holonomy invariant metric on the fibers of the
normal bundle of a Riemannian foliation and, more generally, to
any first order transversally elliptic operator with holonomy
invariant transverse principal symbol. Noncommutative geometry
provides a universal way to develop geometry on $M/\cF$, starting
from the spectral triples associated with this space. To study
such a geometry and investigate its relations with "naive"
geometry of $M/\cF$ (transverse geometry of $\cF$) seems to be a
quite interesting and important problem. Moreover, the language of
noncommutative geometry seems to be very natural and convenient in
the study of many problems of spectral theory and index theory for
differential operators adapted to a foliated structure on a
manifold.

As mentioned above, we will only consider the simplest case of the
leaf space of a Riemannian foliation. Connes and Moscovici in
\cite{Co-M} constructed a spectral triple in a closely related
situation of a compact manifold, equipped with an arbitrary (not
necessarily isometric) action of discrete (pseudo)group. They used
the so-called (transverse) mixed signature operator on the total
space of the (transverse) frame bundle and transversally
hypoelliptic operators. We don't discuss this construction here,
referring the interested reader to \cite{Co-M} (see also
\cite{survey} and references cited therein).

The development of noncommutative geometry of foliations raises
many interesting problems in analysis and geometry in foliated
manifolds. One of our main goals in this paper is to formulate
some of these problems.

Let us describe the contents of the paper. In Section~\ref{s:pdo},
we collect necessary background information on classical
pseudodifferential calculus. In Section~\ref{s:ncspaces}, we
introduce the operator algebras associated with the leaf space
$M/\cF$ of a compact foliated manifold $(M,\cF)$ and with the
cotangent bundle to $M/\cF$.

In Section~\ref{s:transpdo}, we turn to the corresponding quantum
objects associated with the leaf space $M/\cF$. We describe an
appropriate pseudodifferential calculus --- the classes $\Psi
^{m,-\infty}(M,{\mathcal F},E)$ of transversal pseudodifferential
operators on $M$, the corresponding symbolic calculus and their
basic properties. It should be noted that the algebra of symbols
in the transversal pseudodifferential calculus is a noncommutative
algebra. Actually, it is a noncommutative analogue of the algebra
of functions on the cotangent bundle to $M/\cF$, which is
introduced in Section~\ref{s:ncspaces}.

Section~\ref{s:trdyn} is devoted to classical and quantum
dynamical systems on the leaf space $M/\cF$. We introduce
Hamiltonian flows on the cotangent bundle to $M/\cF$ as
one-parameter groups of automorphisms of the associated
noncommutative algebra and formulate the Egorov theorem for
transversally elliptic operators, which provides a relation
between the quantum evolution of transverse pseudodifferential
operators and the corresponding Hamiltonian dynamics on the
cotangent bundle to $M/\cF$ --- the classical evolution of
symbols.

In Section~\ref{s:nc}, we give the definition of a spectral triple
and introduce some geometric objects on the noncommutative space
defined by a spectral triple. We describe spectral triples
associated with the transverse Riemannian geometry of a Riemannian
foliation and give a description of various geometric and analytic
objects determined by these spectral triples in terms of the
classical objects of the transverse geometry of foliations.

We will assume some basic knowledge of foliation theory, referring
the reader to our survey paper \cite{survey} for a summary of
results and, for instance. to the books
\cite{Camacho,Candel-Conlon1,Candel-Conlon2,Godbillon,MMrcun,Molino,M-S,Re,Tondeur}
for different aspects of foliation theory. We also refer the
reader to \cite{survey} and the references cited therein for more
information on noncommutative geometry of foliations.

The author is grateful to N. Azamov and F. Sukochev for very
useful discussions on Dixmier traces and to the referee for very
careful reading and useful remarks.

\section{Preliminaries on pseudodifferential
operators}\label{s:pdo} Pseudodifferential operators are quantum
mechanical observables for a quantum point particle on a compact
manifold. Therefore, they play an important role in our
considerations. For convenience of the reader, we collect in this
Section some necessary facts about pseudodifferential operators
(for more information on pseudodifferential operators see, for
instance, \cite{H3,Taylor,Treves1,Shubin:pdo}).

\subsection{Definition of classes}
Let $U$ be an open subset of $\RR^N$.

\begin{defn}
A function $k\in C^{\infty}(U \times {\mathbb R}^{q},{\mathcal
L}({\mathbb C}^r))$ belongs to the class $S^{m}(U \times {\mathbb
R}^{q}, {\mathcal L}({\mathbb C}^r))$, if, for any multi-indices
$\alpha\in\NN^q $ and $\beta\in\NN^{N}$, there is a constant
$C_{\alpha \beta} > 0$ such that
\[
\|
\partial^{\alpha}_{\eta}
\partial^{\beta}_{x}k(x,\eta )\| \leq C_{\alpha \beta
}(1 +\vert \eta \vert )^{ m-\vert \alpha \vert },\quad x\in
U,\quad \eta \in {\mathbb R}^{q}.
\]
\end{defn}
Here we use notation $|\alpha|=\alpha_1+\alpha_2+\ldots+\alpha_q$
for a multi-index $\alpha\in\NN^q$, and, for a Hilbert space $V$,
${\mathcal L}(V)$ denotes the space of linear bounded maps in $V$.

In the following, we will only consider classical symbols.

\begin{defn}
A function $k\in C^{\infty}(U \times {\mathbb R}^{q}, {\mathcal
L}({\mathbb C}^r))$ is called a classical symbol of order $z\in
{\mathbb C}$, if it can be represented as an asymptotic sum  $$
k(x,\eta)\sim \sum_{j=0}^{\infty} \theta(\eta) k_{z-j}(x,\eta), $$
where $k_{z-j}\in C^{\infty}(U \times ({\mathbb R}^{q}\backslash
\{0\}), {\mathcal L}({\mathbb C}^r))$ are homogeneous in $\eta$ of
degree $z-j$, that is, $$ k_{z-j}(x,t\eta)=t^{z-j}k_{z-j}(x,\eta),
\quad t>0, $$ and $\theta$ is a smooth function in ${\mathbb
R}^{q}$ such that $\theta(\eta)=0$ for $|\eta|\leq 1$,
$\theta(\eta)=1$ for $|\eta|\geq 2$.
\end{defn}

In this definition, the asymptotic equivalence $\sim$ means that,
for any natural $K$,
\[
k - \sum_{j=0}^{K} \theta k_{z-j} \in S^{{\rm Re}\,z-K-1}(U \times
{\mathbb R}^{q}, {\mathcal L}({\mathbb C}^r)).
\]

Consider the $n$-dimensional cube $I^{n}=(0,1)^n$. A classical
symbol $k \in S ^{m} (I ^{n} \times {\mathbb R}^{n}, {\mathcal
L}({\mathbb C}^r))$ defines an operator $A: C^{\infty}_c(I^n,
{\mathbb C}^r) \rightarrow C^{\infty}(I^n, {\mathbb C}^r)$ as
\begin{equation}
\label{pdo0:loc} Au(x)=(2\pi)^{-n} \int e^{i(x-x')\eta}k(x,\eta)
u(x') \,dx'\,d\eta,
\end{equation}
where \(u \in C^{\infty}_{c}(I^{n}, {\mathbb C}^r), x \in I^{n}\).
Denote by $\Psi^{m}(I^n,{\mathbb C}^r)$ the class of operators of
the form (\ref{pdo0:loc}) with $k \in S ^{m} (I ^{n} \times
{\mathbb R}^{n}, {\mathcal L}({\mathbb C}^r))$ such that its
Schwartz kernel is compactly supported in $I^n\times I^n$.

Now let $M$ be a compact $n$-dimensional manifold and $E$ a
complex vector bundle of rank $r$ on $M$. Consider two coordinate
charts on $M$, $\phi: U\to I^n$ and $\phi': U'\to I^n$, endowed
with trivializations of $E$. An operator $A\in\Psi^{m}
(I^n,{\mathbb C}^r)$ determines an operator $A':
C^{\infty}_c(U,\left.E\right|_U) \rightarrow
C^{\infty}_c(U',\left.E\right|_{U'})$, which can be extended in a
trivial way to an operator in $C^{\infty}(M,E)$. The operator
obtained in such a way will be called an elementary operator of
class $\Psi^{m} (M,E)$.

Denote by $\Psi ^{-\infty}(M,E)$ the class of smoothing operators
in $C^{\infty}(M,E)$, i.e., operators with a smooth Schwartz
kernel.

\begin{defn}
The class $\Psi ^{m}(M,E)$ consists of operators $A$, acting in
$C^{\infty}(M,E)$, which can be represented in the form
\[
A=\sum_{i=1}^k A_i + K,
\]
where $A_i$ are elementary operators of
class $\Psi ^{m}(M,E)$, corresponding to pairs $U_i,U'_i$ of
coordinate charts, and $K\in \Psi ^{-\infty}(M,E)$.
\end{defn}

This definition is equivalent to usual definitions of
pseudodifferential operators, but it is more convenient for our
purposes. To see this equivalence, take any finite cover of $M$ by
coordinate charts, $M=\cup_{i=1}^d U_i$. Let $\phi_i\in
C^{\infty}(M), i=1,\ldots,d$ be a partition of unity subordinate
to this cover, $\supp \phi_i\subset U_i,$ and let $\psi_i\in
C^{\infty}(M)$ be such that $\supp \psi_i\subset U_i$ and
$\psi_i\equiv 1$ on $\supp \phi_i$. Then an operator $A\in \Psi
^{m}(M,E)$ is written as
\[
A=\sum_{i=1}^d\psi_iA\phi_i+K, \quad K\in \Psi ^{-\infty}(M,E),
\]
and, for any $i$, $\psi_iA\phi_i$ is an elementary operator of
class $\Psi ^{m}(M,E)$, corresponding to the pair $U_i, U_i$ of
coordinate charts.

A similar definition was used by A. Connes in \cite{Co79} (see
also \cite{M-S}) to introduce the classes of leafwise
pseudodifferential operators on a foliated manifold.

\subsection{Symbolic calculus}
The principal symbol $\sigma_A$ of an elementary operator
$A\in\Psi^{m}(I^n,{\mathbb C}^r)$ of the form (\ref{pdo0:loc}) is
defined to be a smooth matrix-valued function $\sigma_A$ on
$I^n\times ({\mathbb R}^n \backslash \{0\})$ given by
\begin{equation}
\label{princ0} \sigma_A(x,\eta) = k_m(x,\eta),\quad (x,\eta) \in
I^n\times ({\mathbb R}^n \backslash \{0\}),
\end{equation}
where $k_m$ is the homogeneous of degree $m$ component of $k$.

Now let $M$ be a compact $n$-dimensional manifold and $E$ a
complex vector bundle on $M$. Denote by $\pi^*E$ the lift of $E$
to the punctured cotangent bundle $\tilde{T}^*M=T^*M\setminus 0$
under the bundle map $\pi:\tilde{T}^*M\to M$.

The space of all $s\in C^{\infty}(\tilde{T}^*M, {\mathcal
L}(\pi^*E))$, homogeneous of degree $m$ with respect to the
$\RR_+$-multiplication in the fibers of the bundle
$\pi:\tilde{T}^*M\to M$, is denoted by $S^{m}(\tilde{T}^*M,
{\mathcal L}(\pi^*E))$. The linear space
\[
S^{*}(\tilde{T}^*M, {\mathcal L}(\pi^*E))=\bigcup_{m\in\ZZ}
S^{m}(\tilde{T}^*M, {\mathcal L}(\pi^*E))
\]
has the structure of an involutive algebra given by the pointwise
multiplication and the pointwise transposition.

For an operator $A\in \Psi^{m} (M,E)$, the functions defined by
(\ref{princ0}) in any coordinate chart determine a well-defined
element $\sigma_A$ of $S^{m} (\tilde{T}^*M, {\mathcal L}(\pi^*E))$
--- the principal symbol of $A$.

\begin{prop}
The space
\[
\Psi^{*} (M,E)= \bigcup_{m\in\ZZ} \Psi^{m} (M,E)
\]
has the structure of an involutive algebra given by the
composition and transposition of operators. The principal symbol
map
\[
\sigma: \Psi^{*} (M,E)\rightarrow S^{*}(\tilde{T}^*M, {\mathcal
L}(\pi^*E))
\]
is a $\ast$-homomorphism of involutive algebras. In other words:

(1) If $A\in \Psi^{m_1} (M,E)$ and $B\in \Psi^{m_2} (M,E)$,
  then  \(C = AB\) belongs to \(\Psi ^{m_{1}+m_{2}}(M,E)\) and
$\sigma_{AB}=\sigma_A\sigma_B$.

(2) If $A\in \Psi ^{m} (M,E)$, then \(A^{*} \in \Psi ^{m}(M,E)\)
and $\sigma_{A^*}= (\sigma_A)^*$.
\end{prop}

Any $A\in \Psi^{0} (M,E)$ defines a bounded operator in the
Hilbert space $L^2(M,E)$. If $A\in \Psi^{m} (M,E)$ for some $m<0$,
then $A$ is a compact operator in $L^2(M,E)$. Denote by
$\bar{\Psi}^{0}(M,E)$ the closure of $\Psi^{0}(M,E)$ in the
uniform topology of ${\mathcal L}(L^2(M,E))$.

Observe that the algebra $S^{0}(\tilde{T}^*M, {\mathcal
L}(\pi^*E))$ is naturally isomorphic to $C^\infty(S^*M, {\mathcal
L}(\pi^*E))$ and its closure in the uniform topology is isomorphic
to $C(S^*M, {\mathcal L}(\pi^*E))$.

\begin{prop}\label{p:sigma}
$($1$)$ The principal symbol map  $\sigma $ extends by continuity
to a surjective homomorphism
\[
\bar{\sigma}:\bar{\Psi}^{0}(M,E) \rightarrow C(S^*M, {\mathcal
L}(\pi^*E)).
\]

$($2$)$ The ideal $\Ker \bar{\sigma}$ coincides with the ideal
$\cK(L^2(M,E))$ of compact operators in $L^2(M,E)$.
\end{prop}

By Proposition~\ref{p:sigma}, we have a short exact sequence
\[
0 \longrightarrow \cK(L^2(M,E)) \longrightarrow \bar{\Psi}^0(M,E)
\longrightarrow C(S^*M, {\mathcal L}(\pi^*E)) \longrightarrow 0,
\]
which describes the structure of the $C^*$-algebra
$\bar{\Psi}^0(M,E)$ and provides a description of the cosphere
bundle from the operator data
\begin{equation}\label{e:cosphere}
C(S^*M, {\mathcal L}(\pi^*E)) \cong
\bar{\Psi}^0(M,E)/\cK(L^2(M,E)).
\end{equation}

\subsection{The residue trace and zeta-functions} Let $M$ be a compact manifold, $E$
a vector bundle on $M$ and $P\in\Psi^{*}(M,E)$. The residue trace
$\tau(P)$ introduced by Wodzicki \cite{Wo} and Guillemin
\cite{Gu85} is defined as follows. First, the residue form
$\rho_P$ of $P$ is defined in local coordinates as $$ \rho_P =
\left(\int_{|\xi|=1}\Tr p_{-n}(x,\xi)\,d\xi\right) |dx|, $$ where
$p_{-n}(x,\xi)$ is the homogeneous of degree $-n$ ($n=\dim M$) in
$\xi$ component of the complete symbol of $P$. The density
$\rho_P$ turns out to be independent of the choice of a local
coordinate system and, therefore, determines a well-defined
density on $M$. The integral of $\rho_P$ over $M$ is, by
definition, the residue trace $\tau(P)$ of $P$:
\begin{equation}\label{e:wodz}
\tau(P)=(2\pi)^{-n}\int_M\rho_P =(2\pi)^{-n}\int_{S^*M}\Tr
p_{-n}(x,\xi)\, dx d\xi.
\end{equation}
Wodzicki \cite{Wo} showed that $\tau$ is a unique trace on the
algebra $\Psi^*(M,E)$.

Recall that an operator $A\in \Psi^m(M,E)$ is elliptic, if its
principal symbol $\sigma _{A}(x,\xi)$ is invertible for any
$(x,\xi)\in \widetilde{T}^*M$. Examples of elliptic operators are
given by the signature operator $D=d+d^*$ and the Laplace
operator $\Delta=D^2=dd^*+d^*d$ on differential forms on a
compact Riemannian manifold and by the Dirac operator on a
compact Riemannian spin manifold.

\begin{theorem}
Let $A\in\Psi^{m}(M,E)$ be a positive self-adjoint elliptic
operator with the positive definite principal symbol. For any
$Q\in \Psi^{l}(M,E)$, $l\in {\mathbb Z}$, the function $z\mapsto
\tr(QA^{-z})$ is holomorphic for ${\rm Re}\, z> (l+n)/m$ and
admits a (unique) meromorphic extension to ${\mathbb C}$ with at
most simple poles at $z_k=k/m$ with integer $k\leq l+n$. Its
residue at the point $z=z_k$ equals
\[
\underset{z=z_k}{\res} \tr (QA^{-z})=n \tau(QA^{-k/m}).
\]
\end{theorem}

As a consequence, we get the Weyl asymptotic formula for the
eigenvalue distribution function $N(\lambda)$ of a self-adjoint
elliptic operator $A\in\Psi^{m}(M)$ with the positive principal
symbol $\sigma_A$. Let $\lambda_1\leq \lambda_2 \leq \ldots ,
\lambda_m\to +\infty$ be the eigenvalues of $A$ (counted with
multiplicities) and let $\phi_j\in C^\infty(M)$ be a corresponding
orthonormal system of eigenfunctions such that
$A\phi_j=\lambda_j\phi_j, j=1,2,\ldots$. As $\lambda\to +\infty$,
one has
\begin{multline}\label{e:Weyl}
N(\lambda)=\sharp\{j: \lambda_j\leq\lambda\}\\
=(2\pi)^{-n}\lambda^{n/m}{\rm vol}\,\{(x,\xi)\in T^*M:
\sigma_A(x,\xi)\leq 1\}+ O(\lambda^{(n-1)/m}).
\end{multline}

More generally, we have the local Weyl asymptotic formula, which
asserts that, for any $Q\in \Psi^{l}(M)$, $l\in {\mathbb Z}$, one
has as $\lambda\to +\infty$
\begin{multline*}
\sum_{j: \lambda_j\leq\lambda} (Q\phi_j,\phi_j) \\
=(2\pi)^{-n}\lambda^{(l+n)/m} \int_{\{(x,\xi)\in T^*M:
\sigma_A(x,\xi)\leq 1\}}\sigma_Q(x,\xi)\,dx\,d\xi +
O(\lambda^{(l+n-1)/m}).
\end{multline*}

\subsection{Egorov's theorem}
Recall that a classical dynamical system on a compact manifold $M$
(the configuration space) is given by a Hamiltonian flow $f_t$ on
the cotangent bundle $T^*M$ (the phase space) associated with a
classical Hamiltonian $H\in C^\infty(T^*M)$. A quantum dynamical
system on $M$ is given by a one-parameter group of
$\ast$-automorphisms of the algebra $\cL(L^2(M))$:
\[
A\in \cL(L^2(M)) \mapsto A(t)=e^{itP}Ae^{-itP}\in \cL(L^2(M)),
\]
associated with a quantum Hamiltonian $P$, which is a self-adjoint
(unbounded) linear operator in $L^2(M)$. If $P\in \Psi^1(M)$ is a
positive self-adjoint operator and $p$ is its principal symbol,
then the Hamiltonian flow $f_t$ on $T^*M$ associated with $p$ is
called the bicharacteristic flow of $P$. In the case
$P=\sqrt{\Delta_g}\in \Psi^1(M)$, where $\Delta_g$ is the
Laplacian of a Riemannian metric $g$ on $M$, the bicharacteristic
flow of $P$ is the geodesic flow on $T^*M$ associated with $g$.

The Egorov theorem \cite{egorov} relates the quantum evolution of
pseudodifferential operators with the classical dynamics of
principal symbols.

\begin{thm}\label{t:cl-egorov}
Let $M$ be a compact manifold, $E$ a vector bundle on $M$ and
$P\in \Psi^1(M,E)$ a positive self-adjoint pseudodifferential
operator with the positive principal symbol $p$.
\medskip\par
(1) If $A\in \Psi^0(M,E)$, then $A(t)=e^{itP}Ae^{-itP}\in
\Psi^0(M,E)$.
\medskip\par
(2) Moreover, if $E$ is the trivial line bundle and $a\in
S^0(\tilde{T}^*M)$ is the principal symbol of $A$, then the
principal symbol  $a_t\in S^0(\tilde{T}^*M)$ of $A(t)$ is given by
\[
a_t(x,\xi)=a(f_t(x,\xi)), \quad (x,\xi)\in \tilde{T}^*M,
\]
where  $f_t$ is the bicharacteristic flow of $P$.
\end{thm}

What we have described above is a so-called homogeneous
quantization. Its non-homogeneous version, which is associated
with $T^*M$ rather than with $S^*M$, involves Planck's constant
$\hbar$, $\hbar$-dependent pseudodifferential operators and
semiclassical analysis. The corresponding  semiclassical version
of Egorov's theorem is proved in \cite{Robert}.

\section{Some noncommutative spaces associated with the leaf
space}\label{s:ncspaces} In this Section, we will briefly describe
the noncommutative algebras associated with the leaf space of a
foliation. For a more detailed information on various concepts and
facts of noncommutative geometry of foliations, we refer the
reader to a survey \cite{survey} and the bibliography cited
therein.

First, we define a ``nice'' algebra, consisting of functions, on
which all basic operations of analysis are defined. Depending on a
problem in question, one can complete this algebra and obtain a
noncommutative analogue of an appropriate function algebra, for
instance, a von Neumann algebra, an analogue of the algebra of
measurable functions, or a $C^*$-algebra, an analogue of the
algebra of continuous functions, or a smooth algebra, an analogue
of the algebra of smooth functions. The role of a ``nice'' algebra
is played by the algebra $C^\infty_c(G)$ of smooth compactly
supported functions on the holonomy groupoid $G$ of the foliation.
Therefore, we start with the notion of holonomy groupoid of a
foliation.

\subsection{The holonomy groupoid of a foliation}
First, recall the general definition of a groupoid.

\begin{defn}
We say that a set $G$ has the structure of a groupoid with the set
of units $G^{(0)}$, if there are defined maps
\begin{itemize}
\item $\Delta : G^{(0)}\rightarrow G$ (the diagonal map or the unit map);
\item an involution $i:G\rightarrow G$ called the inversion and
written as $i(\gamma)=\gamma^{-1}$;
\item a range map $r:G\rightarrow G^{(0)}$ and a source map $s:G\rightarrow G^{(0)}$;
\item an associative multiplication $m: (\gamma,\gamma')\rightarrow
\gamma\gamma'$ defined on the set
\[
G^{(2)}=\{(\gamma,\gamma')\in G\times G : r(\gamma')=s(\gamma)\},
\]
\end{itemize}
satisfying the conditions
\begin{itemize}
\item $r(\Delta(x))=s(\Delta(x))=x$ and $\gamma\Delta(s(\gamma))=\gamma$,
$\Delta(r(\gamma))\gamma=\gamma$;
\item $r(\gamma^{-1})=s(\gamma)$ and $\gamma\gamma^{-1}=\Delta(r(\gamma))$.
\end{itemize}
\end{defn}

Alternatively, one can define a groupoid as a small category,
where each morphism is an isomorphism.

It is convenient to think of an element $\gamma\in G$ as an arrow
$\gamma :x \to y$, going from $x=s(\gamma)$ to $y=r(\gamma)$.

We will use the standard notation (for $x,y\in G^{(0)}$):
\begin{itemize}
\item $G^x=\{\gamma\in G:r(\gamma)=x\} =r^{-1}(x)$,
\item $G_x=\{\gamma\in G:s(\gamma)=x\} =s^{-1}(x)$,
\item $G_x^y=\{\gamma\in G : s(\gamma)=x, r(\gamma)=y\}$.
\end{itemize}

The holonomy groupoid $G$ of a foliated manifold $(M,{\mathcal
F})$ is defined in the following way. Let $\sim_h$ be an
equivalence relation on the set of continuous leafwise paths
$\gamma:[0,1]\rightarrow M$, setting $\gamma_1\sim_h \gamma_2$, if
$\gamma_1$ and $\gamma_2$ have the same initial and final points
and the same holonomy maps: $h_{\gamma_1} = h_{\gamma_2}$. The
holonomy groupoid $G$ is the set of $\sim_h$-equivalence classes
of leafwise paths. The set of units $G^{(0)}$ is the manifold $M$.
The multiplication in $G$ is given by the product of paths. The
corresponding source and range maps $s,r:G\rightarrow M$ are given
by $s(\gamma)=\gamma(0)$ and $r(\gamma)=\gamma(1)$. Finally, the
diagonal map $\Delta:M\rightarrow G$ takes any $x\in M$ to the
element in $G$ given by the constant path $\gamma(t)=x, t\in
[0,1]$. To simplify the notation, we will identify $x\in M$ with
$\Delta(x)\in G$.

For any $x\in M$ the map $s$ maps $G^x$ on the leaf $L_x$ through
$x$. The group $G^x_x$ coincides with the holonomy group of $L_x$.
The map $s:G^x\rightarrow L_x$ is the covering map associated with
the group $G^x_x$, called the holonomy covering.

The holonomy groupoid $G$ has the structure of a smooth (in
general, non-Hausdorff and non-paracompact) manifold of dimension
$2p+q$. In the following, we will always assume that $G$ is a
Hausdorff manifold.

There is a foliation ${\mathcal G}$ of dimension $2p$ on the
holonomy groupoid $G$. The leaf of $\cG$ through $\gamma\in G$
consists of all $\gamma'\in G$ such that $r(\gamma)$ and
$r(\gamma')$ lie on the same leaf of $\cF$.

\subsection{The noncommutative leaf space of a foliation}
Here we give the intrinsic definition of the operator algebra
associated with a foliated manifold, which uses no additional
choices. It will use the language of half-densities. Indeed, we
will usually consider operators, acting on half-densities, because
their use makes our considerations more natural and simple.

We recall some basic facts concerning densities and integration of
densities (cf., for instance, \cite{Cannas-Wein,geom:asymp}).

\begin{defn}
Let $L$ be an $n$-dimensional linear space and $\cB(L)$ the set of
bases in $L$. An $\alpha$-density on $L$ ($\alpha\in \RR$) is a
function $\rho : \cB(L)\to \CC$ such that, for any $A=(A_{ij})\in
\mathop{GL}(n,\CC)$ and $e=(e_1, e_2,\ldots, e_n)\in \cB(L)$,
\[
\rho(e\cdot A)=|\det A|^\alpha \rho(e),
\]
where $(e\cdot A)_i=\sum_{j=1}^ne_jA_{ji}, i=1,2,\ldots,n$.
\end{defn}

We will denote by $|L|^\alpha$ the space of all $\alpha$-densities
on $L$. For any vector bundle $V$ on $M$, denote by $|V|^{\alpha}$
the associated bundle of $\alpha$-densities, $|V|=|V|^1$.

For any smooth, compactly supported density $\rho$ on a smooth
manifold $M$ there is a well-defined integral $\int_M\rho$,
independent of the fact if $M$ is orientable or not. This fact
allows to define a Hilbert space $L^2(M)$, canonically associated
with $M$, which consists of square integrable half-densities on
$M$. The diffeomorphism group of $M$ acts on $L^2(M)$ by unitary
transformations.

Let $(M,\cF)$ be a compact foliated manifold. Consider the vector
bundle of leafwise half-densities $|T{\mathcal F}|^{1/2}$ on $M$.
Pull back $|T{\mathcal F}|^{1/2}$ to the vector bundles
$s^*(|T{\mathcal F}|^{1/2})$ and $r^*(|T{\mathcal F}|^{1/2})$ on
the holonomy groupoid $G$, using the source map $s$ and the range
map $r$. Define a vector bundle $|T\cG|^{1/2}$ on $G$ as $$
|T\cG|^{1/2}=r^*(|T{\mathcal F}|^{1/2})\otimes s^*(|T{\mathcal
F}|^{1/2}). $$ The bundle $|T\cG|^{1/2}$ is naturally identified
with the bundle of leafwise half-densities on the foliated
manifold $(G,{\mathcal G})$.

The structure of an involutive algebra on
$C^{\infty}_c(G,|T\cG|^{1/2})$ is defined as
\begin{equation}\label{e:algebra}
\begin{aligned}
\sigma_1\ast \sigma_2(\gamma)&=\int_{\gamma_1\gamma_2=\gamma}
\sigma_1(\gamma_1)\sigma_2(\gamma_2),\quad \gamma\in G,\\
\sigma^*(\gamma)&=\overline{\sigma(\gamma^{-1})},\quad \gamma\in
G,
\end{aligned}
\end{equation}
where $\sigma, \sigma_1, \sigma_2\in C^{\infty}_c (G,
|T\cG|^{1/2})$. The formula for $\sigma_1\ast \sigma_2$ should be
interpreted in the following way. If we write $\gamma : x\to y,
\gamma_1: z\to y $ and $\gamma_2: x\to z$, then
\begin{align*}
\sigma_1(\gamma_1)\sigma_2(\gamma_2) \in & |T_y{\mathcal
F}|^{1/2}\otimes |T_z{\mathcal F}|^{1/2}\otimes |T_z{\mathcal
F}|^{1/2}\otimes |T_x{\mathcal F}|^{1/2} \\
 & \cong |T_y{\mathcal
F}|^{1/2}\otimes |T_z{\mathcal F}|^{1}\otimes |T_x{\mathcal
F}|^{1/2},
\end{align*}
and, integrating the $|T_z{\mathcal F}|^{1}$-component
$\sigma_1(\gamma_1)\sigma_2(\gamma_2)$ with respect to $z\in M$,
we get a well-defined section of the bundle $r^*(|T{\mathcal
F}|^{1/2})\otimes s^*(|T{\mathcal F}|^{1/2})=|T\cG|^{1/2}. $

As mentioned above, the algebra $C^{\infty}_c(G,|T\cG|^{1/2})$
plays a role of noncommutative analogue of algebra of functions
on the leaf space $M/\cF$. As we will explain later, this algebra
consists of smooth functions on the leaf space $M/\cF$ in the
sense of noncommutative geometry.

We will also need an analogue of a vector bundle on the leaf space
$M/\cF$ given by a holonomy equivariant vector bundle on $M$. The
corresponding noncommutative analogue of a vector bundle on
$M/\cF$ is given by an appropriate bimodule over
$C^{\infty}_c(G,|T\cG|^{1/2})$, but we don't need this notion
here.

\begin{defn}
A vector bundle $E$ on a foliated manifold $(M,\cF)$ is called
holonomy equivariant, if there is given a representation $T$ of
the holonomy groupoid $G$ of the foliation $\cF$ in the fibers of
$E$, that is, for any $\gamma\in G, \gamma:x\rightarrow y$, there
is defined a linear operator $T(\gamma):E_x\rightarrow E_y$ such
that $T(\gamma_1\gamma_2)=T(\gamma_1)T(\gamma_2)$ for any
$\gamma_1,\gamma_2 \in G$ with $r(\gamma_2)=s(\gamma_1)$.

A Hermitian vector bundle $E$ on a foliated manifold $(M,\cF)$ is
called holonomy equivariant, if it is a holonomy equivariant
vector bundle and the representation $T$ is unitary:
$T(\gamma^{-1})=T(\gamma)^*$ for any $\gamma\in G$.
\end{defn}

Let $E$ be a holonomy equivariant Hermitian vector bundle on a
compact foliated manifold $(M,\cF)$. Any $\sigma\in
C^{\infty}_c(G, |T\cG|^{1/2})$ defines a bounded operator
$R_E(\sigma)$ in the space $C^{\infty}(M,E\otimes |TM|^{1/2})$ of
smooth half-densities on $M$ with values in $E$. For any $u\in
C^{\infty}(M, E\otimes |TM|^{1/2})$, the element $R_E(\sigma)u$ of
$C^{\infty}(M, E\otimes |TM|^{1/2})$ is given by
\[
R_E(\sigma)u(x)=\int_{G^x} \sigma(\gamma)(T\otimes dh^*)(\gamma)
s^*u(\gamma),\quad x\in M,
\]
where $dh^*: s^*(|TM/T\cF|^{1/2}) \to r^*(|TM/T\cF|^{1/2})$ is
induced by the linear holonomy map.

This formula should be interpreted as follows. First, note that
$|TM|^{1/2}\cong |T\cF|^{1/2}\otimes |TM/T\cF|^{1/2}$. We have
\[
s^*u\in C^{\infty}_c(G, s^*(E\otimes |TM|^{1/2}))
\]
and, hence,
\[
\sigma (T\otimes dh^*)\cdot s^*u \in C^{\infty}_c(G, r^*E \otimes
r^*(|T{\mathcal F}|^{1/2})\otimes r^*(|TM/T\cF|^{1/2}) \otimes
s^*(|T{\mathcal F}|)).
\]
The integration of the component in $s^*(|T{\mathcal F}|)$ over
$G^x$, i.e. with a fixed $r(\gamma)=x\in M$, gives a well-defined
section $R_E(\sigma)u$ of $E\otimes |TM|^{1/2}$ on $M$. The
correspondence $\sigma\mapsto R_E(\sigma)$ defines a
representation $R_E$ of the algebra $C^{\infty}_c(G,
|T\cG|^{1/2})$ in $C^{\infty}(M,E\otimes |TM|^{1/2})$.

\subsection{The noncommutative cotangent bundle to the leaf
space}\label{s:cot} Like in classical theory, the cotangent bundle
to the leaf space of a foliation and its quantization will play a
very important role in our considerations. In this section, we
describe the corresponding noncommutative object. To do this, we
will follow the construction of the cotangent bundle $T^*B$ to the
base $B$ from the cotangent bundle $T^*M$ to the total space $M$
for a fibration $M\to B$ (as explained in \cite{egorgeo}, this
construction can be considered as a particular case of the
foliation reduction in symplectic geometry) and, when it will be
necessary, switch to noncommutative algebras.

Assume that the foliation $\cF$ is Riemannian. Let $N^*{\mathcal
F}=\{\nu\in T^*M:\langle\nu,X\rangle=0\ \mbox{\rm for any}\ X\in
T\cF\}$ denote the conormal bundle to ${\mathcal F}$. If $(x,y)\in
I^p\times I^q$ denotes the local coordinates in a foliated chart
$\phi: U\to I^p \times I^q$ and $(x,y,\xi,\eta) \in I^p \times I^q
\times {\RR}^p\times {\RR}^q$ the local coordinates in the
corresponding chart on $T^*M$, then the subset $N^*{\mathcal
F}\cap \pi^{-1}(U)=U_1$ (here $\pi: T^*M\to M$ is the bundle map)
is given by $\xi=0$.

There is the natural lift of $\cF$ to a foliation ${\mathcal F}_N$
on $N^*{\mathcal F}$ called the horizontal (or linearized)
foliation. The coordinate chart $\phi_n: N^*{\mathcal F}\to I^p
\times I^q\times {\RR}^q$ determined by a foliated coordinate
chart $\phi$ on $M$ is a foliated chart for ${\mathcal F}_N$ with
plaques given by the level sets $y={\rm const}, \eta={\rm const}$.

The leaf $\widetilde{L}_{\nu}$ of ${\mathcal F}_N$ through a point
$\nu\in N^*{\mathcal F}$ consists of all points of the form
$dh_{\gamma}^{*}(\nu)$ with $\gamma\in G$ such that
$r(\gamma)=\pi(\nu)$. It is diffeomorphic to the holonomy covering
$G^x$ of the leaf $L_x, x=\pi(\nu)$ of ${\mathcal F}$ through $x$.
Each leaf of the linearized foliation ${\mathcal F}_N$ has trivial
holonomy.

The leaf space $N^*{\mathcal F}/{\mathcal F}_N$ of the foliation
${\mathcal F}_N$ can be considered as the cotangent bundle to the
leaf space $M/\cF$. This holds in the case when the foliation is
given by a fibration, but, in general, the leaf space is
singular, and we will consider the associated operator algebras.

The holonomy groupoid $G_{{\mathcal F}_N}$ of the foliation
${\mathcal F}_N$ is described as follows: $$ G_{{\mathcal F}_N}=
\{(\gamma,\nu)\in G\times N^*{\mathcal F} : r(\gamma)=\pi(\nu)\}
$$ with the source map $s_N:G_{{\mathcal F}_N}\rightarrow
N^{*}{\mathcal F}, s_N(\gamma,\nu)=dh_{\gamma}^{*}(\nu)$, the
range map $r_N:G_{{\mathcal F}_N}\rightarrow N^{*}{\mathcal F},
r_N(\gamma,\nu)=\nu$ and the composition
$(\gamma,\nu)(\gamma',\nu') = (\gamma\gamma',\nu)$ defined in the
case when $\nu'=dh_{\gamma}^{*}(\nu)$. The projection $\pi:
N^{*}{\mathcal F}\rightarrow M$ induces a map $\pi_G:G_{{\mathcal
F}_N}\rightarrow G$ by the formula $ \pi_G(\gamma,\nu)=\gamma,
(\gamma,\nu)\in G_{{\mathcal F}_N}. $ Denote by ${\mathcal G}_N$
the natural foliation on $G_{{\mathcal F}_N}$.

Taking into account the fact that $\tilde{N}^*\cF$ is noncompact,
we introduce the space $C^{\infty}_{prop}(G_{{\mathcal F}_N},
|T{\mathcal G}_N|^{1/2})$, which consists of all properly
supported elements $k\in C^{\infty}(G_{{\mathcal F}_N},
|T{\mathcal G}_N|^{1/2})$ (this means that the restriction of
$r:G_{\cF_N}\to \tilde{N}^*\cF$ to $\supp k$ is a proper map).
Then one can introduce the structure of involutive algebra on
$C^{\infty}_{prop}(G_{{\mathcal F}_N}, |T{\mathcal G}_N|^{1/2})$,
using the formulas (\ref{e:algebra}).

The algebra $C^{\infty}_{prop}(G_{{\mathcal F}_N}, |T{\mathcal
G}_N|^{1/2})$ plays a role of a noncommutative analogue of algebra
of functions on the cotangent bundle to the leaf space $M/\cF$.

\section{Transverse pseudodifferential calculus}\label{s:transpdo}
Now we turn to the quantum objects associated with the leaf space
of a compact foliated manifold $(M,\cF)$. We need an appropriate
pseudodifferential calculus, the classes $\Psi
^{m,-\infty}(M,{\mathcal F},E)$ of transversal pseudodifferential
operators, which was developed in \cite{noncom}. In this section,
we recall the definition of classes $\Psi ^{m,-\infty}(M,{\mathcal
F},E)$ and their basic properties. These classes can be considered
as a slight generalization of the algebra of Fourier integral
operators associated to a coisotropic submanifold of a symplectic
manifold \cite{GS79} in the particular case when the symplectic
manifold is $T^*M$ and the coisotropic submanifold is the conormal
bundle $N^*\cF$ to $\cF$.

\subsection{Definition of classes}
Consider the $n$-di\-men\-si\-o\-nal cube $I^{n} = I^{p} \times
I^{q}$ equipped with a trivial foliation, whose leaves are
$I^{p}\times \{y\}$, $y\in I^{q}$. The coordinates in $I^{n}$ will
be denoted by $(x,y)$, $x\in I^{p}$, $y \in I^{q}$, and the dual
coordinates by $(\xi ,\eta )$, $\xi \in {\RR}^{p}$, $\eta \in
{\RR}^{q}$.

A classical symbol $k \in S ^{m} (I ^{p} \times I^p\times
I^q\times {\mathbb R}^{q}, {\mathcal L}({\mathbb C}^r))$ defines
an operator
\[
A: C^{\infty}_c(I^n, {\mathbb C}^r) \rightarrow C^{\infty}(I^n,
{\mathbb C}^r)
\]
as
\begin{equation}
\label{pdo:loc} Au(x,y)=(2\pi)^{-q} \int
e^{i(y-y')\eta}k(x,x',y,\eta) u(x',y') \,dx'\,dy'\,d\eta,
\end{equation}
where \(u \in C^{\infty}_{c}(I^{n}, {\mathbb C}^r), x \in I^{p}, y
\in I^{q}\). Denote by $\Psi^{m,-\infty}(I^n,I^p,{\mathbb C}^r)$
the class of operators of the form (\ref{pdo:loc}) with $k \in S
^{m} (I ^{p} \times I^p\times I^q\times {\mathbb R}^{q}, {\mathcal
L}({\mathbb C}^r))$ such that its Schwartz kernel is compactly
supported in $I^n\times I^n$.

Let $(M,\cF)$ be a compact foliated manifold, $\dim M=n$, $\dim
\cF=p$, $p+q=n$, and let $E$ be a vector bundle of rank $r$ on
$M$. Let $\phi: U\to I^p\times I^q, \phi': U'\to I^p\times I^q$ be
two foliated charts, $\pi=pr_{nq}\circ\phi: U\to \RR^q$,
$\pi'=pr_{nq}\circ\phi': U'\to \RR^q$ the corresponding
distinguished maps. The foliated charts $\phi$, $\phi'$ are called
compatible, if, for any $m\in U$ and $m'\in U'$ with
$\pi(m)=\pi'(m')$, there is a leafwise path $\gamma$ from $m$ to
$m'$ such that the corresponding holonomy map $h_{\gamma}$ takes
the germ $\pi_m$ of $\pi$ at $m$ to the germ $\pi'_{m'}$ of $\pi'$
at $m'$.

If $\phi: U \subset M\to I^p\times I^q, \phi': U' \subset M \to
I^p\times I^q$ are compatible foliated charts on $M$ endowed with
trivializations of $E$, then an operator
$A\in\Psi^{m,-\infty}(I^n,I^p,{\mathbb C}^r)$ defines an operator
$A': C^{\infty}_c(U,\left.E\right|_U)\rightarrow
C^{\infty}_c(U',\left.E\right|_{U'})$, which can be extended in a
trivial way to an operator in $C^{\infty}(M,E)$. The operator
obtained in such a way will be called an elementary operator of
class $\Psi^{m,-\infty} (M,{\mathcal F},E)$.

\begin{defn}
The class $\Psi ^{m,-\infty}(M,{\mathcal F},E)$ consists of
operators $A$, acting in $C^{\infty}(M,E)$, which can be
represented in the form
\[
A=\sum_{i=1}^k A_i + K,
\]
where $A_i$ are elementary operators of
class $\Psi ^{m,-\infty}(M,{\mathcal F},E)$, corresponding to
pairs $\phi_i,\phi'_i$ of compatible foliated charts, and $K\in
\Psi ^{-\infty}(M,E)$.
\end{defn}

\subsection{Symbolic calculus}
The principal symbol $\sigma_A$ of an elementary operator
$A\in\Psi^{m,-\infty}(I^n,I^p,{\mathbb C}^r)$ given by
(\ref{pdo:loc}) is defined to be the matrix-valued half-density
$\sigma_A$ on $I^p\times I^p\times I^q\times ({\mathbb R}^q
\backslash \{0\})$ given by
\begin{multline}
\label{princ} \sigma_A(x,x',y,\eta) = k_m(x,x',y,\eta)
|dx\,dx'|^{1/2},\\ (x,x',y,\eta) \in I^p\times I^p\times I^q\times
({\mathbb R}^q \backslash \{0\}),
\end{multline}
where $k_m$ is the homogeneous of degree $m$ component of $k$.

Let $(M,\cF)$ be a compact foliated manifold and let $E$ be a
Hermitian vector bundle on $M$. Denote by $\pi^*E$ the lift of $E$
to the punctured conormal bundle $\tilde{N}^*\cF=N^*\cF\setminus
0$ under the map $\pi:\tilde{N}^*\cF\to M$. Denote by ${\mathcal
L}(\pi^*E)$ the vector bundle on $G_{{\mathcal F}_N}$, whose fiber
at a point $(\gamma,\nu)\in G_{{\mathcal F}_N}$ consists of all
linear maps from $(\pi^*E)_{s_N(\gamma,\nu)}$ to
$(\pi^*E)_{r_N(\gamma,\nu)}$, where, for any $\nu \in
\tilde{N}^*\cF$, $(\pi^*E)_{\nu}$ denotes the fiber of $\pi^*E$
at $\nu$. One can introduce the structure of involutive algebra
on the space $C^{\infty}_{prop}(G_{{\mathcal F}_N}, {\mathcal
L}(\pi^*E)\otimes |T{\mathcal G}_N|^{1/2})$ of all properly
supported sections of the vector bundle ${\mathcal
L}(\pi^*E)\otimes |T{\mathcal G}_N|^{1/2}$ on $G_{{\mathcal
F}_N}$ by formulas similar to (\ref{e:algebra}).

The space of all sections $s\in C^{\infty}_{prop}(G_{{\mathcal
F}_N},{\mathcal L}(\pi^*E)\otimes |T{\mathcal G}_N|^{1/2})$,
homogeneous of degree $m$ with respect to the
$\RR_+$-multiplication in the fibers of the bundle
$\pi:\tilde{N}^*\cF\to M$, is denoted by $S^{m}(G_{{\mathcal
F}_N},{\mathcal L}(\pi^*E)\otimes |T{\mathcal G}_N|^{1/2})$. The
space $$ S^*(G_{{\mathcal F}_N},{\mathcal L}(\pi^*E)\otimes
|T{\mathcal G}_N|^{1/2})=\bigcup_{m\in\ZZ} S^m(G_{{\mathcal
F}_N},{\mathcal L}(\pi^*E)\otimes |T{\mathcal G}_N|^{1/2}) $$ is a
subalgebra of $C^{\infty}_{prop}(G_{{\mathcal F}_N}, {\mathcal
L}(\pi^*E)\otimes |T{\mathcal G}_N|^{1/2})$.

Let $\phi: U\subset M\to I^p\times I^q, \phi': U' \subset M \to
I^p\times I^q$ be two compatible foliated charts on $M$ endowed
with trivializations of $E$. Then the corresponding coordinate
charts $\phi_n: U_1\subset N^*\cF\to I^p\times I^q \times \RR^q,
\phi_n': U_1^\prime \subset N^*\cF \to I^p\times I^q \times \RR^q$
are compatible foliated charts on the foliated manifold
$(N^*\cF,\cF_N)$ endowed with obvious trivializations of $\pi^*E$.
Thus, there is a foliated chart $\Gamma_N: W(\phi_n,\phi'_n)
\subset G_{{\mathcal F}_N}\to I^p\times I^p\times I^q\times
{\mathbb R}^q$ on the foliated manifold $(G_{\cF_N},{\cG}_N)$.

For $A\in \Psi^{m,-\infty} (M, {\mathcal F},E)$, the
half-densities defined by (\ref{princ}) in any foliated chart
$W(\phi_n,\phi'_n)$ determine a well-defined element $\sigma_A$
of $S^m(G_{{\mathcal F}_N},{\mathcal L}(\pi^*E)\otimes |T{\mathcal
G}_N|^{1/2})$ --- the principal symbol of $A$.

\begin{prop}
\label{prop} The space
\[
\Psi^{*,-\infty}(M,{\mathcal F},E)= \bigcup_{m\in\ZZ}
\Psi^{m,-\infty}(M,{\mathcal F},E)
\]
has the structure of an involutive algebra given by the
composition and transposition of operators. The principal symbol
map
\[
\sigma: \Psi^{*,-\infty}(M,{\mathcal F},E)\rightarrow
S^*(G_{{\mathcal F}_N},{\mathcal L}(\pi^*E)\otimes |T{\mathcal
G}_N|^{1/2})
\]
is a $\ast$-homomorphism of involutive algebras.
\end{prop}

Recall that the principal symbol of a pseudodifferential operator
$P$ acting in $C^\infty(M,E)$ is a well-defined section of the
bundle $\cL(\pi^*E)$ on $\widetilde{T}^*M$, where $\pi : T^*M\to
M$ is the natural projection.

\begin{defn}
The transversal principal symbol $\sigma _{P}$ of an operator
$P\in \Psi ^{m}(M,E)$ is the restriction of its principal symbol
$p_m$ to $\widetilde{N}^{\ast}{\mathcal F}$.
\end{defn}

\begin{prop}
\label{module} If $A\in  \Psi ^{\mu}(M,E)$ and $B\in  \Psi
^{m,-\infty}(M,{\mathcal F},E)$, then $AB$ and $BA$ belong to
$\Psi ^{\mu+m,-\infty}(M,{\mathcal F},E)$ and
\begin{align*}
\sigma_{AB}(\gamma,\nu)&=\sigma_A(\nu)\sigma_B(\gamma,\nu), \quad
(\gamma,\nu)\in G_{{\mathcal F}_N},\\[6pt]
\sigma_{BA}(\gamma,\nu)&=
\sigma_B(\gamma,\nu)\sigma_A(dh_{\gamma}^*(\nu)), \quad
(\gamma,\nu)\in G_{{\mathcal F}_N}.
\end{align*}
\end{prop}

Suppose that $E$ is holonomy equivariant, that is, there is an
action $T(\gamma):E_x\rightarrow E_y, \gamma \in G,
\gamma:x\rightarrow y$ of the holonomy groupoid $G$ in the fibers
of $E$. Then the bundle ${\mathcal L}(\pi^*E)$ on $N^*\cF$ is
holonomy equivariant with the corresponding action
$\operatorname{ad} T$ of the holonomy groupoid $G_{\cF_N}$ in the
fibers of ${\mathcal L}(\pi^*E)$.

\begin{defn}\label{d:hinv}
The transversal principal symbol $\sigma_P$ of an operator
$P\in\Psi^m(M,E)$ is holonomy invariant, if, for any leafwise path
$\gamma $ from $x$ to $y$ and for any $\nu \in N^{
\ast}_{y}{\mathcal F}$, the following identity holds: $$ {\rm
ad}\,T(\gamma,\nu)[\sigma _{P}(dh_{\gamma}^{\ast}(\nu ))] = \sigma
_{P}(\nu ). $$
\end{defn}

The assumption of the existence of a positive order
pseudodifferential operator with a holonomy invariant transversal
principal symbol on a foliated manifold imposes sufficiently
strong restrictions on geometry of the foliation. An example of an
operator with a holonomy invariant transverse principal symbol is
given by the transverse signature operator on a Riemannian
foliation.

There is a canonical embedding $$ i:C^{\infty}_{prop}(G_{{\mathcal
F}_N},|T{\mathcal G}_N|^{1/2})\hookrightarrow
C^{\infty}_{prop}(G_{{\mathcal F}_N},{\mathcal L}(\pi^*E)\otimes
|T{\mathcal G}_N|^{1/2}),$$ which takes any $k\in
C^{\infty}_{prop}(G_{{\mathcal F}_N},|T{\mathcal G}_N|^{1/2})$ to
$i(k)=k \pi^*T$. We will identify $C^{\infty}_{prop}(G_{{\mathcal
F}_N},|T{\mathcal G}_N|^{1/2})$ with its image in
$C^{\infty}_{prop}(G_{{\mathcal F}_N},{\mathcal L}(\pi^*E)\otimes
|T{\mathcal G}_N|^{1/2})$ under the map $i$.

\begin{defn}
An operator $P\in \Psi^{m,-\infty}(M,{\mathcal F},E)$ is said to
have a scalar principal symbol, if its principal symbol belongs to
$C^{\infty}_{prop}(G_{{\mathcal F}_N},|T{\mathcal G}_N|^{1/2})$.
\end{defn}

Denote by $\Psi_{sc}^{m,-\infty}(M,{\mathcal F},E)$ the set of all
operators $K\in \Psi^{m,-\infty}(M,{\mathcal F},E)$ with a scalar
principal symbol. Observe that, for any $k\in
C^{\infty}_c(G,|T{\cG}|^{1/2})$, the operator $R_E(k)$ belongs to
$\Psi^{0,-\infty}_{sc}(M,\cF,E)$, and
\[
\sigma(R_E(k))=\pi_G^*k\in C^{\infty}_{prop}(G_{{\mathcal
F}_N},|T{\mathcal G}_N|^{1/2}).
\]

\begin{prop}
Let $(M,\cF)$ be a compact foliated manifold and $E$ a holonomy
equivariant vector bundle. If
$A\in\Psi_{sc}^{m,-\infty}(M,{\mathcal F},E)$, and $P\in\Psi
^{\mu}(M,E)$ has a holonomy invariant transversal principal
symbol, then $[A,P]$ is in $\Psi^{m+\mu-1,-\infty}(M,{\mathcal
F})$.
\end{prop}

Any $A\in \Psi^{0,-\infty} (M,{\mathcal F},E)$ defines a bounded
operator in the Hilbert space $L^2(M,E)$. Denote by
$\bar{\Psi}^{0,-\infty}(M,{\mathcal F},E)$ the closure of
$\Psi^{0,-\infty}(M,{\mathcal F},E)$ in the uniform topology of
${\mathcal L}(L^2(M,E))$.

For any $\nu\in \tilde{N}^*\cF$, there is a natural
$\ast$-representation $R_\nu$ of the algebra $S^{0}(G_{{\mathcal
F}_N}, {\mathcal L}(\pi^*E)\otimes |T{\mathcal G}_N|^{1/2})$ in
$L^2(G^\nu_{{\cF}_N}, s_N^*(\pi^*E))$. Thus, for any $k \in
S^{0}(G_{{\mathcal F}_N}, {\mathcal L}(\pi^*E)\otimes |T{\mathcal
G}_N|^{1/2})$, the continuous operator family
\[
\{R_\nu(k)\in \cL(L^2(G^\nu_{{\cF}_N}, s_N^*(\pi^*E))):\nu\in
\tilde{N}^*\cF\}
\]
defines a bounded operator in $L^2(G_{{\cF}_N},s_N^*(\pi^*E))$. We
will identify $k$ with the corresponding bounded operator in
$L^2(G_{{\cF}_N},s_N^*(\pi^*E))$ and denote by
$\bar{S}^{0}(G_{{\mathcal F}_N}, {\mathcal L}(\pi^*E)\otimes
|T{\mathcal G}_N|^{1/2})$ the closure of $S^{0}(G_{{\mathcal
F}_N}, {\mathcal L}(\pi^*E)\otimes |T{\mathcal G}_N|^{1/2})$ in
the uniform topology of $\cL(L^2(G_{{\cF}_N},s_N^*(\pi^*E)))$.

\begin{prop}[\cite{egorgeo}]
\label{seq} $($1$)$ The symbol map $$ \sigma :
\Psi^{0,-\infty}(M,{\mathcal F},E)\rightarrow S^{0}(G_{{\mathcal
F}_N},{\mathcal L}(\pi^*E)\otimes |T{\mathcal G}_N|^{1/2}) $$
extends by continuity to a homomorphism $$
\bar{\sigma}:\bar{\Psi}^{0,-\infty}(M,{\mathcal F},E) \rightarrow
\bar{S}^{0}(G_{{\mathcal F}_N},{\mathcal L}(\pi^*E)\otimes
|T{\mathcal G}_N|^{1/2}). $$

$($2$)$ The ideal $\Ker \bar{\sigma}$ contains the ideal of
compact operators in $L^2(M,E)$.
\end{prop}

We have much less information on the principal symbol map in
transverse pseudodifferential calculus. For instance, answers to
the following questions are unknown.

\begin{quest}\label{q:seq1}
Is the principal symbol map $\bar{\sigma}$ surjective?
\end{quest}

\begin{quest}\label{q:seq2}
Under which conditions is the principal symbol map $\bar{\sigma}$
injective?
\end{quest}

Let us make some comments. Recall that the representation $R_E$
determines an inclusion
\[
C^{\infty}_c(G,|T{\cG}|^{1/2})\longrightarrow
\Psi^{0,-\infty}_{sc}(M,\cF,E)
\]
and the restriction of $\sigma$ to $C^{\infty}_c(G,
|T{\cG}|^{1/2})$ is the identity map, if we identify
$C^{\infty}_c(G,|T{\cG}|^{1/2})$ with its image in
$C^{\infty}_{prop}(G_{{\mathcal F}_N},|T{\mathcal G}_N|^{1/2})$ by
the map $\pi_G^*$ induced by the projection $\pi_G :G_{\cF_N}\to
G$. Passing to the completions, we will get a homomorphism
\[
\pi_E : C^{\ast}_{E}(G) \rightarrow C^{\ast}_{r}(G)
\]
where $C^{*}_E(G)$ is the closure of
$R_E(C^{\infty}_c(G),|T{\cG}|^{1/2})$ in the uniform operator
topology of ${\cL}(L^2(M,E))$ and $C^{\ast}_{r}(G)$ is the reduced
$C^*$-algebra of $G$. By \cite{F-Sk}, this homomorphism is
surjective, but, in general, is not injective. It is injective for
any $E$ if the groupoid $G$ is amenable (cf., for instance,
\cite{F-Sk} and also \cite{claire-jean}). Therefore, if $G$ is not
amenable, we cannot expect that $\bar\sigma$ is injective.

\subsection{The residue trace and zeta-functions}
There is an analogue of the Wodzicki-Guillemin residue trace for
operators from $\Psi^{m,-\infty}(M,{\mathcal F},E)$ \cite{noncom},
which is defined as follows. First, note that it suffices to
define the residue trace for elementary operators of class
$\Psi^{m,-\infty}(M,{\mathcal F},E)$. For
$P\in\Psi^{m,-\infty}(I^n,I^p,{\mathbb C}^r)$, define the residue
form $\rho_P$ as $$ \rho_P = \left(\int_{|\eta|=1} \Tr
k_{-q}(x,x,y,\eta)\,d\eta\right) |dx dy|, $$ and the residue trace
$\tau(P)$ as $$ \tau(P)=(2\pi)^{-q}\int_{|\eta|=1}\Tr
k_{-q}(x,x,y,\eta)\, dx dy d\eta, $$ where $k_{-q}$ is the
homogeneous of degree $-q$ component of the complete symbol $k$ of
$P$.

For any $P\in \Psi^{m,-\infty}(M,{\mathcal F},E)$, its residue
form $\rho_P$ is a well-defined density on $M$, and the residue
trace $\tau(P)$ is obtained by the integration of $\rho_P$ over
$M$: $$ \tau(P)=(2\pi)^{-q} \int_{M}\rho_P. $$

\begin{defn}
A pseudodifferential operator $P\in \Psi^m(M,E)$ is called
transversally elliptic, if its transversal principal symbol
$\sigma _{P}(\nu)$ is invertible for any $\nu\in
\widetilde{N}^*{\mathcal F}$.
\end{defn}

\begin{theorem}[\cite{noncom}]\label{t:zeta}
Let $A\in\Psi^{m}(M,E)$ be a transversally elliptic operator with
a positive transversal principal symbol. Suppose that the operator
$A$, considered as an unbounded operator in the Hilbert space
$L^2(M,E)$, is essentially self-adjoint on the initial domain
$C^{\infty}(M,E)$, and its closure is an invertible and positive
operator.

For any $Q\in \Psi^{l,-\infty}(M,{\mathcal F},E)$, $l\in {\mathbb
Z}$, the function $z\mapsto \tr(QA^{-z})$ is holomorphic for ${\rm
Re}\, z> l+q/m$ and admits a (unique) meromorphic extension to
${\mathbb C}$ with at most simple poles at $z_k=k/m$ with integer
$k\leq l+q$. Its residue at the point $z=z_k$ equals
\[
\underset{z=z_k}{\res} \tr (QA^{-z})=q \tau(QA^{-k/m}).
\]
\end{theorem}

One can easily derive from Theorem~\ref{t:zeta} a Weyl type
asymptotic formula for the distributional spectrum distribution
function of a positive transversally elliptic operator with a
positive transversal principal symbol, as well as an asymptotic
expansion for its distributional heat trace.

\begin{prob}
To extend Theorem~\ref{t:zeta} to the case when the symbol of
$Q\in \bar{\Psi}^{0,-\infty}(M,\cF,E)$ (which belongs to
$C^{\infty}(G_{{\mathcal F}_N},{\mathcal L}(\pi^*E)\otimes
|T{\mathcal G}_N|^{1/2})$) is not properly supported.
\end{prob}

One can expect that this result holds in the case when the symbol
of $Q$ is exponentially decreasing at infinity. If this is the
case, this fact can be considered as a sort of quantum ergodic
theorem for foliations, and the rate of the exponential decay
could be related with a version of (tangential) entropy of
foliations.

\subsection{Adiabatic limits and noncommutative Weyl formula}\label{s:adiab} Let
$(M,{\mathcal F})$ be a closed foliated manifold, $\dim M = n$,
$\dim {\mathcal F} = p$, $p+q=n$, endowed with a Riemannian metric
$g_M$. Then we have a decomposition of the tangent bundle to $M$
into a direct sum $TM=F\oplus H$, where $F=T{\mathcal F}$ is the
tangent bundle to $\cF$ and $H=F^{\bot}$ is the orthogonal
complement of $F$, and the corresponding decomposition of the
metric: $ g_{M}=g_{F}+g_{H}$. Define a one-parameter family
$g_{h}$ of Riemannian metrics on $M$ by
\begin{equation}\label{e:gh}
g_{h}=g_{F} + {h}^{-2}g_{H}, \quad 0 < h \leq 1.
\end{equation}

For any $h>0$, consider the Laplace operator $\Delta_{h}$ on
differential forms defined by the metric $g_h$. It is a
self-adjoint, elliptic, differential operator with the positive,
scalar principal symbol in the Hilbert space $L^2(M,\Lambda
T^{*}M,g_h)$ of square integrable differential forms on $M$,
endowed with the inner product induced by $g_h$, which has
discrete spectrum. In \cite{adiab}, the asymptotic behavior of the
trace of $f(\Delta_h)$ when $h\to 0$ was studied for any $f\in
S(\RR)$. Such asymptotic limits are called adiabatic limits after
Witten.

It turns out that this asymptotic spectral problem can be
considered as a semiclassical spectral problem for a Schr\"odinger
operator on the leaf space $M/\cF$, and the resulting asymptotic
formula for the trace of $f(\Delta_h)$ can be written in the form
of the semiclassical Weyl formula for a Schr\"odinger operator on
a compact Riemannian manifold, if we replace the classical objects
entering to this formula by their noncommutative analogues.

To demonstrate this, first, transfer the operators $\Delta_{h}$ to
the fixed Hilbert space $L^2(M,\Lambda T^{*}M)=L^{2}(M,\Lambda
T^{*}M, g_M)$, using an isomorphism $\Theta_h$ from
$L^{2}(M,\Lambda T^{*}M , g_{h})$ to $L^{2}(M,\Lambda T^{*}M)$
defined as follows. With respect to a bigrading on $\Lambda
T^{*}M$ given by
\[
\Lambda^k T^{*}M=\bigoplus_{i=0}^{k}\Lambda^{i,k-i}T^{*}M, \quad
\Lambda^{i,j}T^{*}M=\Lambda^{i}F^{*}\otimes \Lambda^{j}H^{*},
\]
we have
\[
\Theta_{h}u = h^{j}u, \quad u \in L^{2}(M,\Lambda^{i,j}T^{*}M ,
g_{h}).
\]
The operator $\Delta_h$ in $L^{2}(M,\Lambda T^{*}M , g_{h})$
corresponds under the isometry $\Theta_{h}$ to the operator
$L_{h}= \Theta_{h}\Delta_h\Theta_{h}^{-1}$ in $L^{2}(M,\Lambda
T^{*}M)$.

With respect to the bigrading of $\Lambda T^*M$, the de Rham
differential $d$ can be written as
\[
d=d_F+d_H+\theta,
\]
where
\begin{enumerate}
\item  $ d_F=d_{0,1}: C^{\infty}(M,\Lambda^{i,j}T^*M)\to
C^{\infty}(M,\Lambda^{i,j+1}T^*M)$ is the tangential de Rham
differential, which is a first order tangentially elliptic
operator, independent of the choice of $g_M$;
\item $d_H=d_{1,0}: C^{\infty}(M,\Lambda^{i,j}T^*M)\to
C^{\infty}(M,\Lambda^{i+1,j}T^*M)$ is the transversal de Rham
differential, which is a first order transversally elliptic
operator;
\item $\theta=d_{2,-1}: C^{\infty}(M,\Lambda^{i,j}T^*M)\to
C^{\infty}(M,\Lambda^{i+2,j-1}T^*M)$ is a zero order differential
operator.
\end{enumerate}

In the case when ${\mathcal F}$ is a Riemannian foliation and
$g_{M}$ is a bundle-like metric, one can show that the leading
term in the asymptotic expansion of the trace of $f(\Delta_h)$ or,
that is the same, of the trace of $f(L_h)$ as $h\to 0$ coincides
with the leading term in the asymptotic expansion of the trace of
$f(\bar{L}_h)$ as $h\to 0$, where
\[
\bar{L}_h=\Delta_F + h^2\Delta_H,
\]
$\Delta_F=d_Fd^*_F+d^*_Fd_F$ is the tangential Laplacian and
$\Delta_H=d_Hd^*_H+d^*_Hd_H$ is the transverse Laplacian.

Now observe that the operator $\bar{L}_h$ has the form of a
Schr\"odinger operator on the leaf space $M/\cF$, where $\Delta_H$
plays a role of the Laplace operator, and $\Delta_F$ a role of the
operator-valued potential on $M/\cF$.

Recall that in the case of a Schr\"odinger operator $H_h$ on a
compact Riemannian manifold $X$ with a matrix-valued potential
$V\in C^\infty(X,{\mathcal L}(E))$, where $E$ is a
finite-dimensional Euclidean space and $V(x)^{*}=V(x)$:
\[
H_h=-h^2\Delta +V(x),\quad x\in X,
\]
the corresponding asymptotic formula (the semiclassical Weyl
formula) has the following form:
\[ \operatorname{tr}
f(H_h)=(2\pi)^{-n}h^{-n}\int_{T^*X} \operatorname{Tr}
f(p(x,\xi))\,dx\,d\xi+o(h^{-n}),\quad h\rightarrow 0+,
\]
where $p\in C^\infty(T^*X,\cL(E))$ is the principal $h$-symbol of
$H_h$:
\begin{equation*}
p(x,\xi)=|\xi|^2+V(x),\quad (x,\xi)\in T^{*}X.
\end{equation*}
Now let us show how the asymptotic formula for the trace of
$f(\Delta_h)$ in the adiabatic limit can be written in a similar
form, using noncommutative geometry. First, we define the
principal $h$-symbol of $\Delta_h$. Denote by $g_N$ the Riemannian
metric on $N^*{\mathcal F}$ induced by the Riemannian metric on
$M$. The principal $h$-symbol of $\Delta_h$ is a tangentially
elliptic operator in $C^{\infty}(N^*{\mathcal F},\pi^{*} \Lambda
T^{*}M)$ given by
\[
\sigma_h(\Delta_h) = \Delta_{{\mathcal F}_N}+g_N,
\]
where $\Delta_{{\mathcal F}_N}$ is the lift of the tangential
Laplacian $\Delta_F$ to a tangentially elliptic (relative to
${\mathcal F}_N$) operator $\Delta_{{\mathcal F}_N}$ in
$C^{\infty}(N^*{\mathcal F},\pi^{*} \Lambda T^{*}M)$, and $g_N$
denotes the multiplication operator by $g_N \in C^\infty(N^*\cF)$.
(Observe that $g_N$ coincides with the transversal principal
symbol of $\Delta_H$.)

We will consider $\sigma_h(\Delta_h)$ as a family of elliptic
operators along the leaves of the foliation ${\mathcal F}_N$. For
any function $f\in C^\infty_c({\mathbb R})$, the operator
$f(\sigma_h(\Delta_h))$ belongs to the twisted foliation
$C^*$-algebra $C^*(N^*\cF, {\mathcal F}_N,\pi^{*}\Lambda T^{*}M)$,
which is the noncommutative analogue of continuous differential
forms on the leaf space $N^*{\cF}/\cF_N$, the cotangent bundle to
$M/\cF$.

Then we replace the usual integration over $T^{*}X$ and the matrix
trace $\operatorname{Tr}$ by the integration in the sense of the
noncommutative integration theory given by the trace
$\operatorname{tr}_{{\mathcal F}_N}$ on the twisted foliation
$C^*$-algebra, which is defined by the canonical transverse
Liouville measure for the symplectic foliation ${\mathcal F}_N$.
One can show that the value of this trace on
$f(\sigma_h(\Delta_h))$ is finite.

\begin{theorem}[\cite{adiab}]
\label{ad:main} For any $f\in C^\infty_c({\mathbb R})$, the
asymptotic formula holds:
\begin{equation}\label{e:adiab}
\operatorname{tr} f(\Delta_{h}) =(2\pi)^{-q}h^{-q}
\operatorname{tr}_{{\mathcal F}_N} f(\sigma_h(\Delta_h))
+o(h^{-q}),\quad h\rightarrow 0.
\end{equation}
\end{theorem}

Observe that the formula (\ref{e:adiab}) makes sense for an
arbitrary, not necessarily Riemannian, foliation. Therefore, it is
quite reasonable to conjecture that it holds in such generality.

\begin{conj}
Let $\cF$ be an arbitrary foliation on a compact Riemannian
manifold. In the above notation, for any function $f\in
C^\infty_c({\mathbb R})$, the asymptotic formula holds:
\[
\operatorname{tr} f(\Delta_F + h^2\Delta_H) =(2\pi)^{-q}h^{-q}
\operatorname{tr}_{{\mathcal F}_N} f(\Delta_{{\mathcal F}_N}+g_N)
+o(h^{-q}),\quad h\rightarrow 0.
\]
\end{conj}

To extend the above conjecture to the Laplace operator on $M$, we
can try to use the corresponding signature operators.

\begin{conj}
Let $\cF$ be an arbitrary foliation on a compact Riemannian
manifold. For any even function $f\in C^\infty_c({\mathbb R})$,
the asymptotic formula holds:
\[
\operatorname{tr} f(D_F + h D_H) =(2\pi)^{-q}h^{-q}
\operatorname{tr}_{{\mathcal F}_N} f(D_{{\mathcal
F}_N}+\sigma(D_H)) +o(h^{-q}),\quad h\rightarrow 0,
\]
where $D_F=d_F+d^*_F$ is the leafwise signature operator on $M$,
$D_{{\mathcal F}_N}$ is the correspodning leafwise (relative to
${\mathcal F}_N$) signature operator on $N^*{\cF}$,
$D_H=d_H+d^*_H$ is the transverse signature operator on $M$,
$\sigma(D_H)$ is the transverse principal symbol of $D_H$
(considered as a multiplication operator on $N^*\cF$).
\end{conj}

\section{Transverse dynamics}\label{s:trdyn}
\subsection{Transverse Hamiltonian flows}
In this Section, we will discuss classical dynamical systems on
the leaf space of a foliation. To give their definition, we will
proceed as in Section~\ref{s:cot}. We start with a dynamical
system on the cotangent bundle to the total manifold, satisfying
some symmetry assumptions (like holonomy invariance relative to
the foliation), and try to construct the corresponding dynamical
system on the cotangent bundle to the base. This construction can
be also considered as a particular case of the foliation reduction
in symplectic geometry (see \cite{egorgeo}). Since, in our case,
the base is, in general, a singular object, we pass eventually to
the corresponding operator algebras.

Let $(M,\cF)$ be a compact foliated manifold, and let $p$ be a
homogeneous of degree one function defined in some conic
neighborhood of $\tilde{N}^*{\mathcal F}$ in $\tilde{T}^*M$ such
that its restriction to $\tilde{N}^*{\mathcal F}$ is constant
along the leaves of $\cF_N$. Take any function ${\tilde p}\in
S^1(\tilde{T}^*M)$, which coincides with $p$ in some conic
neighborhood of $\tilde{N}^*{\mathcal F}$. Denote by $X_{\tilde
p}$ the Hamiltonian vector field on $T^*M$ with the Hamiltonian
$\tilde p$. For any $\nu\in \tilde{N}^*{\mathcal F}$, the vector
$X_{\tilde p}(\nu)$ is tangent to $\tilde{N}^*{\mathcal F}$.
Therefore, the Hamiltonian flow $\tilde{f}_t$ with the Hamiltonian
$\tilde{p}$ preserves $\tilde{N}^*{\mathcal F}$. Denote by $f_t$
its restriction to $N^*{\mathcal F}$. One can show that the vector
field $X_{\tilde p}$ on $\tilde{N}^*{\mathcal F}$ is an
infinitesimal transformation of the foliation $\cF_N$, and,
therefore, the flow $f_t$ preserves the foliation $\cF_N$.

It follows from the fact that $X_{\tilde p}$ is an infinitesimal
transformation of $\cF_N$ that there exists a unique vector field
$\cH_p$ on $G_{{\mathcal F}_N}$ such that $ds_N(\cH_p)=X_{\tilde
p}$ and $dr_N(\cH_p)=X_{\tilde p}$. Let $F_t$ be the flow on
$G_{{\mathcal F}_N}$ defined by $\cH_p$. Then $s_N\circ
F_t=f_t\circ s_N$, $r_N\circ F_t=f_t\circ r_N$ and the flow $F_t$
preserves $\cG_N$.

\begin{defn}\label{d:trans}
The transverse Hamiltonian flow of $p$ is the one-parameter group
$F^*_t$ of automorphisms of the involutive algebra
$C^{\infty}_{prop}(G_{{\mathcal F}_N},|T{\mathcal G}_N|^{1/2})$,
induced by the action of $F_t$.
\end{defn}

This definition can be easily seen to be independent of the choice
of $\tilde{p}$.

\subsection{Egorov theorem for transversally elliptic operators}
Let $(M,\cF)$ be a compact foliated manifold, $E$ a Hermitian
vector bundle on $M$ and $D\in\Psi^1(M,E)$ a self-adjoint
transversally elliptic operator in $L^2(M,E)$. Suppose that $D^2$
has the scalar principal symbol and the holonomy invariant
transversal principal symbol. By the spectral theorem, the
operator $\langle D\rangle=(D^2+I)^{1/2}$ defines a strongly
continuous group $e^{it\langle D\rangle}$ of bounded operators in
$L^2(M,E)$. Consider the one-parameter group $\Phi_t$ of
$\ast$-automorphisms of the algebra ${\mathcal L}(L^2(M,E))$
defined as
\begin{displaymath}
\Phi_t(T)=e^{i t\langle D\rangle}Te^{-i t\langle D\rangle}, \quad
T\in {\mathcal L}(L^2(M,E)).
\end{displaymath}

Let $a_2\in S^2(\tilde{T}^*M)$ be the principal symbol of $D^2$
and let $F^*_t$ be the transverse Hamiltonian flow on
$C^{\infty}_{prop}(G_{{\mathcal F}_N},|T{\mathcal G}_N|^{1/2})$
associated with $\sqrt{a_2}$.

Any scalar operator $P\in\Psi^m(M)$, acting on half-densities, has
the subprincipal symbol, which is a globally defined, homogeneous
of degree $m-1$, smooth function on $T^*M\setminus 0$, given in
local coordinates by
\begin{equation}\label{e:subprincipal1}
p_{sub}=p_{m-1}-\frac{1}{2i}\sum_{j=1}^n\frac{\partial^2p_m}{\partial
x_j\partial \xi_j}.
\end{equation}
Note that $p_{sub}=0$, if $P$ is a real self-adjoint differential
operator of even order.

\begin{thm}[\cite{egorgeo}]
\label{Egorov} Let $D\in\Psi^1(M,E)$ be a self-adjoint
transversally elliptic operator in $L^2(M,E)$ such that $D^2$ has
the scalar principal symbol and the holonomy invariant transversal
principal symbol. Let $K\in \Psi^{m,-\infty}(M,{\mathcal F},E)$.
\medskip\par
(1) There is a $K(t)\in\Psi^{m,-\infty}(M,{\mathcal F},E)$ such
that, for any $s$ and $r$, the family $\langle
D\rangle^r(\Phi_t(K)-K(t))\langle D\rangle^{-s}, t\in \RR,$ is a
smooth family of trace class operators in $L^2(M,E)$.
\medskip\par
(2) If, in addition, $E$ is the trivial line bundle, the
subprincipal symbol of $D^2$ vanishes, and $k\in
S^{m}(G_{{\mathcal F}_N},|T{\mathcal G}_N|^{1/2})$ is the
principal symbol of $K$, then the principal symbol $k(t)\in
S^m(G_{{\mathcal F}_N},|T{\mathcal G}_N|^{1/2})$ of $K(t)$ is
given by $k(t)=F^*_t(k)$.
\end{thm}

\begin{prob}\label{p:egorov}
To extend the second statement of Theorem~\ref{Egorov} to the case
when $E$ is an arbitrary vector bundle.
\end{prob}

\subsection{Noncommutative dynamical entropy} In this section we raise a
question, which is very interesting and highly nontrivial even in
the case of compact Riemannian manifold.

So we start with a compact Riemannian manifold $(M,g)$. Recall
that $\bar{\Psi}^0(M)$ denotes the closure of the algebra
$\Psi^0(M)$ in the uniform operator topology in $\cL(L^2(M))$.
Consider the one-parameter group $\Phi_t$ of $\ast$-automorphisms
of the $C^*$-algebra $\bar{\Psi}^0(M)$ defined as
\[
\Phi_t(T)=e^{i t\sqrt{\Delta_g}}Te^{-i t\sqrt{\Delta_g}}, \quad
T\in \bar{\Psi}^0(M),
\]
where $\Delta_g$ is the Laplace operator associated with $g$. Let
$F_t$ denote the geodesic flow on the cosphere bundle $S^*M$ and
$F^*_t$ the induced action on $C(S^*M)$. By the classical Egorov
theorem, Theorem~\ref{t:cl-egorov}, we have the commutative
diagram
\[
  \begin{CD}
\bar{\Psi}^0(M) @>\Phi_t>> \bar{\Psi}^0(M)\\ @V\bar{\sigma}VV
@VV\bar{\sigma}V
\\ C(S^*M) @>F^*_t>> C(S^*M)
  \end{CD}
\]

\begin{prob}
To define a quantum topological entropy $h(\Phi_t)$ of the
noncommutative geodesic flow $\Phi_t$ so that it is related with
the classical topological entropy $h(F_t)$ of the geodesic flow
$F_t$.
\end{prob}

Some very interesting recent results related to this question
were obtained by D. Kerr \cite{Kerr03,Kerr04}.

Now we extend this conjecture to the foliation case.

\begin{prob}
In notation of Theorem~\ref{Egorov}, to define a (classical)
topological entropy $h(F^*)$ of the transverse geodesic flow
$F^*_t$ and a (quantum) topological entropy $h(\Phi)$ of the
noncommutative geodesic flow $\Phi_t$ so that there are relations
between these two notions of entropy.
\end{prob}

\subsection{Noncommutative symplectic geometry}
Based on the ideas of the deformation theory of Gerstenhaber
\cite{Gersten63}, Xu \cite{Xu} and Block and Getzler
\cite{Block-Ge} introduced an analogue of the Poisson bracket in
noncommutative geometry. Namely, they defined a Poisson structure
on an algebra $A$ as a Hochschild $2$-cocycle $P\in Z^2(A,A)$ such
that $P\circ P$ is a Hochschild $3$-coboundary, $P\circ P\in
B^3(A,A)$. In other words, a Poisson structure on $A$ is given by
a linear map $P : A\otimes A\to A$ such that
\begin{multline}\label{e:p1}
(\delta P)(a_1,a_2,a_3)\equiv a_1P(a_2,a_3)-P(a_1a_2,a_3)\\
+P(a_1,a_2a_3)-P(a_1,a_2)a_3=0,
\end{multline}
and there is a $2$-cochain $P_1:A\otimes A\to A$ such that
\begin{multline}\label{e:p2}
P\circ P(a_1,a_2,a_3)\equiv P(a_1,P(a_2,a_3))-P(P(a_1,a_2),a_3)\\
=a_1P_1(a_2,a_3)- P_1(a_1a_2,a_3)+P_1(a_1,a_2a_3)-P_1(a_1,a_2)a_3.
\end{multline}
The identity (\ref{e:p1}) is an analogue of the Jacobi identity
for a Poisson bracket, and the identity (\ref{e:p2}) is an
analogue of the Leibniz rule.

Block and Getzler \cite{Block-Ge} defined a Poisson structure on
the operator algebra $C^\infty_c(G, |T\cG|^{1/2})$ of a
transversally symplectic foliation $\cF$ in the case when the
normal bundle $\tau$ to $\cF$ has a basic connection $\nabla $
(recall that a basic connection on $\tau$ is a holonomy invariant
adapted connection), in particular, when $\cF$ is Riemannian. A
natural example of a transversally symplectic Riemannian foliation
is given by the linearized foliation $\cF_N$ on the conormal
bundle $N^*\cF$ to a {\em Riemannian} foliation $\cF$. So the
construction of Block and Getzler can be applied in this case, and
we get a natural noncommutative Poisson structure on
$C^{\infty}_{prop}(G_{{\mathcal F}_N}, |T{\mathcal G}_N|^{1/2})$.

\begin{prob}
To define the notion of noncommutative Hamiltonian flow on a
noncommutative algebra so that the transverse Hamiltonian flows on
$C^{\infty}_{prop}(G_{{\mathcal F}_N}, |T{\mathcal G}_N|^{1/2})$
would be noncommutative Hamiltonian flows.
\end{prob}

\begin{prob}
To construct (strict) deformation quantization of the algebra
$C^{\infty}_{prop}(G_{{\mathcal F}_N}, |T{\mathcal G}_N|^{1/2})$
(in the sense of Rieffel
\cite{Rieffel88,Rieffel89,Rieffel93,Rieffel94}).
\end{prob}

We refer to \cite{Underhill,Liu-Quian,Landsman93, Carinena} for
some results on quantization of the cotangent bundle and to
\cite{tang} for some recent results on deformation quantization
of symplectic groupoids.

\subsection{Quantum ergodicity} It is well-known that there are
relationships between dynamical properties of the geodesic flow of
a compact Riemannian manifold $(M,g)$ and asymptotic properties of
the eigenvalues and the eigenfunctions of the corresponding
Laplace operator $\Delta_g$. This phenomenon was first discovered
in \cite{Schnirelman74} (see also \cite{Colin85,Zelditch87}).

\begin{thm}[\cite{Schnirelman74}]
Let $(M,g)$ be a compact Riemannian manifold. Let $\lambda_1\leq
\lambda_2\leq\lambda_3\leq\cdots, \lambda_j\to+\infty$ be the
eigenvalues of the associated Laplacian $\Delta_g$ (counted with
multiplicities) and $\varphi_j\in C^\infty(M)$ the corresponding
orthonormal system of eigenfunctions:
\[
  \Delta_g\varphi_j=\lambda_j\varphi_j.
\]
Consider the spectrum distribution function
\[
N(\lambda)=\sharp\{j: \sqrt{\lambda_j}\leq \lambda\}.
\]
If the geodesic flow $G_t$ on $S^*M$ is ergodic, then, for $A\in
\Psi^0(M)$ with the principal symbol $\sigma_A$:
\[
\lim_{\lambda\to+\infty}\frac{1}{N(\lambda)}\sum_{\sqrt{\lambda_j}\leq
\lambda} (A\varphi_j,\varphi_j)=\frac{1}{{\rm
vol}\,(S^*M)}\int_{S^*M}\sigma_A d\mu,
\]
where $d\mu$ is the Liouville measure on $S^*M$.
\end{thm}

The corresponding semiclassical result is due to Helffer, Martinez
and Robert \cite{HMR87}. The development of these results led to
the notions of quantum ergodicity and quantum mixing (see, for
instance, \cite{Sunada97,Zelditch:erg,Zelditch:mix}, and
\cite{Zelditch05} for a recent survey) and belongs to a very
active field of current research in spectral theory of
differential operators and mathematical physics called quantum
chaos.

In Section~\ref{s:adiab}, we have seen that adiabatic limits for
the spectrum of the Laplace operator on a Riemannian foliated
manifold can be naturally considered as semiclassical spectral
problems on the leaf space of the foliation. Therefore, the
following problem is quite natural and its solution would provide
a natural generalization of the results mentioned above to this
setting.

\begin{prob}
To relate dynamical properties of the transverse geodesic flow of
a Riemannian foliation on a compact manifold and asymptotic
properties of the eigenvalues and eigenfunctions of the
corresponding Laplacian in the adiabatic limit.
\end{prob}

\section{Transverse Riemannian geometry}\label{s:nc}
\subsection{Spectral triples}
According to \cite{Co,Co-M,Sp-view}, the initial datum of
noncommutative differential geometry is a spectral triple (or an
unbounded Fredholm module).

\begin{defn}
A spectral triple $({\mathcal A}, {\mathcal H}, D)$ consists of an
involutive algebra ${\mathcal A}$, a Hilbert space ${\mathcal H}$
equipped with a $\ast$-representation of ${\mathcal A}$ (we will
identify an element $a\in \cA$ with the corresponding operator in
$\mathcal H$), and an (unbounded) self-adjoint operator $D$ in
${\mathcal H}$ such that
\medskip\par
1. for any $a\in {\mathcal A}$, the operator $a(D-i)^{-1}$ is a
compact operator in ${\mathcal H}$;

2. for any $a\in {\mathcal A}$, the operator $[D,a]$ is bounded in
$\cH$.
\end{defn}

A spectral triple is supposed to contain the basic geometric
information on Riemannian geometry of the corresponding
geometrical object. In particular, the operator $D$ can be
considered as an analog of Riemannian metric.

We will consider two basic examples of spectral triples:

\subsubsection{Spectral triples associated with compact Riemannian manifolds}
The classical Riemannian geometry is described by the spectral
triple $(\cA,{\mathcal H},D)$ associated with a compact Riemannian
manifold $(M,g)$:
\begin{enumerate}
\item ${\mathcal A}$ is the algebra $C^{\infty}(M)$
of smooth functions on $M$;
\item ${\mathcal H}$ is the space
$L^2(M,\Lambda^*T^*M)$ of differential $L^2$-forms on $M$, on
which the algebra ${\mathcal A}$ acts by multiplication;
\item $D$ is the signature operator $d+d^*$.
\end{enumerate}

\subsubsection{Spectral triples associated with Riemannian
foliations \cite{noncom,egorgeo}} Let $(M,\cF)$ be a compact
foliated manifold. Assume that $\cF$ is Riemannian, and take a
bundle-like metric $g_M$ on $M$. Let $H=F^{\bot}$ be the
orthogonal complement of $F=T{\mathcal F}$ with respect to $g_M$.
Let:
\begin{enumerate}
\item ${\mathcal A}=C^{\infty}_c(G)$;
\item ${\mathcal H}$ is the Hilbert space $L^2(M, \Lambda^{*} H^{*})$
of transverse differential forms;
\item $D$ is the transverse signature operator $d_H+d^*_H$.
\end{enumerate}

More generally, we will consider spectral triples associated with
transversally elliptic operators, acting in sections of a holonomy
equivariant Hermitian vector bundle $E$:
\medskip\par
(T1) ${\mathcal A}=C^{\infty}_c(G)$;

(T2) ${\mathcal H}$ is the Hilbert space $L^2(M,E)$ of $L^2$
sections of $E$ equipped with the action of $\cA$ given by $R_E$;

(T3) $D$ is a first order self-adjoint transversally elliptic
operator, acting in $C^\infty(M,E)$, with the holonomy invariant
transversal principal symbol such that $D^2$ is self-adjoint and
has the scalar principal symbol.

\subsection{Smooth spectral triples}
First, we will describe the noncommutative analogue of a smooth
structure on a topological manifold, the notion of smooth
subalgebra of a $C^*$-algebra, and explain why the operator
algebra $C^\infty_c(G, |T{\cG}|^{1/2})$ associated with a compact
foliated manifold $(M,\cF)$ consists of smooth functions on the
leaf space $M/\cF$ in the noncommutative sense.

Suppose that $A$ is a $C^*$-algebra and $A^+$ is the algebra
obtained by adjoining the unit to $A$. Suppose that $\cA$ is a
$*$-subalgebra of the algebra $A$ and $\cA^+$ is the algebra
obtained by adjoining the unit to $\cA$

\begin{defn}
We say that $\cA$ is a smooth subalgebra of $A$, if:

(1) $\cA$ is a dense $*$-subalgebra of $A$;

(2) $\cA$ is stable under the holomorphic functional calculus,
that is, for any $a\in \cA^+ $ and for any function $f$,
holomorphic in a neighborhood of the spectrum of $a$ (considered
as an element of the algebra $A^+$) $f(a) \in \cA^+$.
\end{defn}

Suppose that ${\cA}$ is a dense $*$-subalgebra of a $C^*$-algebra
$A$, endowed with the structure of a Fr\'echet algebra whose
topology is finer than the topology induced by the topology of
${A}$. By \cite[Lemma 1.2]{Schw}), $\cA$ is a smooth subalgebra of
$A$ if and only if $\cA$ is spectral invariant, that is,
$\cA^+\cap GL(A^+) = GL(\cA^+)$, where $GL(\cA^+)$ and $GL(A^+)$
denote the group of invertibles in $\cA^+$ and $A^+$ respectively.

A spectral triple $({\mathcal A}, {\mathcal H}, D)$ determines a
natural smooth subalgebra in $\cL(\cH)$. Let $\langle
D\rangle=(D^2+I)^{1/2}$. Denote by $\delta$ the (unbounded)
differentiation on ${\mathcal L}({\mathcal H})$ given by
\begin{equation}
\label{derivative} \delta(T)=[\langle D\rangle,T],\quad T\in
\Dom\,\delta\subset {\mathcal L}({\mathcal H}).
\end{equation}
We say that $P\in {\rm OP}^{\alpha}$ if and only if $P\langle
D\rangle^{-\alpha}\in \bigcap_n {\rm Dom}\;\delta^n$. In
particular, ${\rm OP}^{0}=\bigcap_n {\rm Dom}\;\delta^n$. Then
${\rm OP}^0$ is a smooth subalgebra of $\cL(\cH)$ (see, for
instance, \cite[Theorem 1.2]{Ji}).

\begin{defn}\cite{Co-M,Sp-view}
We will say that a spectral triple $({\mathcal A}, {\mathcal H},
D)$ is smooth (or $QC^\infty$ as in \cite{CPRS1}), if, for any
$a\in {\mathcal A}$, we have the inclusions $a, [D,a] \in  {\rm
OP}^{0}$.
\end{defn}

The fact that a spectral triple $({\mathcal A}, {\mathcal H}, D)$
is smooth means that $\cA$ consists of smooth functions on the
corresponding geometric space in the sense of noncommutative
geometry. In particular, for the spectral triple associated with
a compact Riemannian manifold $M$, ${\rm OP}^0\cap C(M)$
coincides with $C^\infty(M)$ (observe that here one can take as
$\cA$ any involutive algebra, which consists of Lipschitz
functions and is dense in $C(M)$).

Let $({\mathcal A}, {\mathcal H}, D)$ be a smooth spectral triple.
Denote by ${\mathcal B}$ the algebra generated by all elements of
the form $\delta^n(a)$, where $a\in{\mathcal A}$ and $n\in
{\mathbb N}$. Thus, $\cB$ is the smallest subalgebra in ${\rm
OP}^{0}$, which contains $\cA$ and is invariant under the action
of $\delta$.

Denote by ${\rm OP}_0^{0}$ the space of all $P\in {\rm OP}^{0}$
such that $\langle D\rangle^{-1}P$ and $P\langle D\rangle^{-1}$
are compact operators in $\cH$. If the algebra $\cA$ has unit,
then ${\rm OP}_0^{0}={\rm OP}^{0}$. By the definition of a
spectral triple, $\cA\subset {\rm OP}^{0}_0$.

\begin{defn}\cite{egorgeo}
We will say that a spectral triple $({\mathcal A}, {\mathcal H},
D)$ is $QC_0^\infty$, if it is smooth and the associated
subalgebra ${\mathcal B}$ is contained in ${\rm OP}_0^{0}$.
\end{defn}

This notion has a natural geometric interpretation. If the algebra
$\cA$ has no unit, we can consider the corresponding
noncommutative space as a noncompact space. The fact that,  for
$a\in \cA$, the operator $a(D-i)^{-1}$ is a compact operator in
${\mathcal H}$ means that $a$ considered as a function on the
corresponding noncommutative space vanishes at infinity. The
condition $\cB\subset {\rm OP}^{0}_0$ means that the elements of
$\cA$ vanish at infinity along with all its derivatives of
arbitrary order.

\begin{thm}\cite{egorgeo}
Any spectral triple defined in (T1), (T2), (T3) is $QC_0^\infty$.
\end{thm}

\subsection{Dimension and dimension spectrum}
As we have been mentioned above, the dimension of a compact
Riemannian manifold can be seen from the Weyl asymptotic formula
for the eigenvalues of the corresponding Laplace (or the
signature) operator (cf. (\ref{e:Weyl})). This fact motivates the
next definition.

For a compact operator $T$ in a Hilbert space $\cH$, denote by
$\mu_1(T)\geq \mu_2(T)\geq \ldots$ the singular numbers of $T$,
that is, the eigenvalues of the operator $|T|=(T^*T)^{1/2}$.
Recall that the Schatten-von Neumann ideal $\cL^{p}(\cH), 1\leq
p<\infty,$ consists of all $T\in \cK(\cH)$ such that
\[
\sum_{n=1}^\infty\mu_n(T)^p<\infty.
\]
The elements of $\cL^{1}(\cH)$ are called trace class operators.
For any $T\in \cL^{1}(\cH)$, its trace is defined as
\[
\tr T=\sum_{n=1}^\infty\mu_n(T).
\]

\begin{defn}\label{d:dim}
A spectral triple $({\mathcal A}, {\mathcal H}, D)$ is called
$p$-summable (or $p$-di\-men\-si\-o\-nal), if, for any $a\in
{\mathcal A}$, the operator $a(D-i)^{-1}$ belongs to
$\cL^p({\mathcal H})$.

A spectral triple $({\mathcal A}, {\mathcal H}, D)$ is called
finite-dimensional, if it is $p$-summable for some $p$.

The greatest lower bound of all $p$'s, for which a
finite-dimensional spectral triple is $p$-summable, is called the
dimension of the spectral triple.
\end{defn}

The spectral triple associated with a compact Riemannian manifold
$(M,g)$ is finite-dimensional, and the dimension of this spectral
triple coincides with the dimension of $M$.

The dimension of spectral triples associated with a Riemannian
foliation $\cF$ is equal to the codimension of $\cF$.

If we are looking at a geometrical space as a union of pieces of
different dimensions, this notion of dimension of the
corresponding spectral triple gives only an upper bound on
dimensions of various pieces. To take into account lower
dimensional pieces of the space under consideration, Connes and
Moscovici \cite{Co-M} suggested that the correct notion of
dimension is given not by a single real number $d$ but by a subset
${\rm Sd}\subset {\CC}$, which is called the dimension spectrum.

\begin{defn}\cite{Co-M,Sp-view}
A spectral triple $({\mathcal A},{\mathcal H},D)$ has the discrete
dimension spectrum ${\rm Sd}\subset{\mathbb C}$, if ${\rm Sd}$ is
a discrete subset in ${\mathbb C}$, the triple is smooth, and, for
any $b\in {\mathcal B}$, the distributional zeta-function
$\zeta_b(z)$ of $\langle D\rangle $ given by $$ \zeta_b(z)=\tr
b\langle D\rangle^{-z}, $$ is defined in the half-plane
$\{z\in\CC:{\rm Re}\, z>d\}$ and extends to a holomorphic function
on ${\mathbb C}\backslash {\rm Sd}$ such that the function
$\Gamma(z)\zeta_b(z)$ is rapidly decreasing on the vertical lines
$z=s+it$ for any $s$ with ${\rm Re}\,s
>0$.

The dimension spectrum is said to be simple, if the singularities
of $\zeta_b(z)$ at $z\in {\rm Sd}$ are at most simple poles.
\end{defn}

The spectral triple associated with a compact Riemannian manifold
has the discrete dimension spectrum, which is contained in $\{v\in
{\mathbb N}:v\leq n=\dim M\}$ and is simple.

\begin{theorem}[\cite{noncom}]
A spectral triple given by (T1), (T2), (T3) has the discrete
dimension spectrum ${\rm Sd}$, which is contained in $\{v\in
{\mathbb N}:v\leq q={\rm codim }\,\cF\}$ and is simple.
\end{theorem}

\subsection{The Dixmier trace and the Riemannian volume form}
In \cite{Dixmier66}, Dixmier introduced a nonstandard trace
$\operatorname{Tr}_\omega$ on the algebra $\cL(\cH)$. Consider the
ideal $\cL^{1+}(\cH)$ in the algebra of compact operators
$\cK(\cH)$, which consists of all $T\in \cK(\cH)$ such that
\[
\sup_{N\in\NN}\frac{1}{\ln N}\sum_{n=1}^N\mu_n(T)<\infty.
\]
For any invariant mean $\omega$ on the amenable group of upper
triangular $2\times 2$-matrices, Dixmier constructed a linear form
$\lim_\omega$ on the space ${\ell}^\infty(\NN)$ of bounded
sequences, which coincides with the limit functional $\lim$ on the
subspace of convergent sequences. The trace
$\operatorname{Tr}_\omega$ is defined for a positive operator
$T\in \cL^{1+}(\cH)$ as
\[
\operatorname{Tr}_\omega(T)=\lim_\omega \frac{1}{\ln N}
\sum_{n=1}^N\mu_n(T).
\]
This trace is non-normal and vanishes on the trace class
operators.

Let $M$ be a compact manifold and $E$ a vector bundle on $M$. As
shown in \cite{Co-action} (cf. also \cite{Gracia:book}), any
operator $P\in\Psi^{-n}(M,E)$ ($n=\dim M$) belongs to the ideal
$\cL^{1+}(L^2(M,E))$, the Dixmier trace
$\operatorname{Tr}_\omega(P)$ does not depend on the choice of
$\omega$ and coincides with the value of the residue trace
$\tau(P)$: for any invariant mean $\omega$,
\[
\operatorname{Tr}_\omega(P)=\tau(P).
\]

For the spectral triple $(\cA,\cH,\cD)$ associated with a compact
Riemannian manifold $(M,g)$, the above results imply the formula
\begin{equation}\label{e:Dix}
\int_Mf\,dx=c(n)\operatorname{Tr}_\omega(f|D|^{-n}), \quad f\in
\cA,
\end{equation}
where $c(n)=2^{(n-[n/2])} \pi^{n/2} \Gamma(\frac{n}{2}+1)$ and
$dx$ denotes the Riemannian volume form on $M$. Thus, the Dixmier
trace $\operatorname{Tr}_\omega$ can be considered as a proper
noncommutative generalization of the integral.

A similar relation of the Dixmier trace $\operatorname{Tr}_\omega$
with the transverse Riemannian volume form associated with a
Riemannian foliation relies on the following conjecture, which
precise formulation have been clarified after our discussions with
N. Azamov and F. Sukochev.

\begin{conj}
Let $(M,\cF)$ be a compact foliated manifold and $E$ a vector
bundle on $M$. Any $P\in \Psi^{-q,-\infty}(M,{\mathcal F},E)$
($q={\rm codim}\,\cF$) belongs to $\cL^{1+}(L^2(M,E))$, the
Dixmier trace $\operatorname{Tr}_\omega(P)$ does not depend on the
choice of $\omega$ and coincides with the value of the residue
trace $\tau(P)$.
\end{conj}

From the other side, if we will consider the residue trace $\tau$
instead of the Dixmier trace $\operatorname{Tr}_\omega$ as the
noncommutative integral, we get the following analog of the
formula (\ref{e:Dix}).

\begin{prop}
Let $(\cA,\cH,D)$ be the spectral triple associated with a
Riemannian foliation $(M,\cF)$. For any $k\in \cA$, we have
\begin{equation}\label{e:int}
\tau(R_E(k)\langle D\rangle^{-q}) =\frac{q}{\Gamma(\frac{q}{2}+1)}
\int_Mk(x)\,dx.
\end{equation}
\end{prop}
Here $k(x)\,dx$ means the product of the restriction of $k$ to
$M$, which is a leafwise density on $M$, and the transverse volume
form of $\cF$. Observe that the right hand side of (\ref{e:int})
coincides (up to some multiple) with the value of the von Neumann
trace $\tr_{\cF}$ given by the transverse Riemannian volume of
$\cF$ due to the noncommutative integration theory \cite{Co79}:
\[
\tr_{\cF}(k)=\int_M k(x)\,dx, \quad k\in C^{\infty}_c
(G,|T{\cG}|^{1/2}).
\]
Recall that $C^{*}_E(G)$ denotes the closure of
$R_E(C^{\infty}_c(G,|T{\cG}|^{1/2}))$ in the uniform operator
topology of ${\cL}(L^2(M,E))$, and $\pi_E : C^{\ast}_{E}(G)
\rightarrow C^{\ast}_{r}(G)$ is the natural projection. A
remarkable observation related with the formula (\ref{e:int}) is
that its right hand side as a functional on $C^{\ast}_{E}(G)$
depends only on $\pi_E(k)$. In particular, for any $k\in \ker
\pi_E $, we have
\[
\underset{z=-q}{\res} \tr (R_E(k)\langle D \rangle^{-z})
=\tau(R_E(k)\langle D \rangle^{-q}) = 0.
\]
One can interpret this fact in the following way. Let us think of
an involutive ideal ${\cI}$ in $C^{*}_E(G)$ as a subset of our
spectrally defined geometrical space. Then if ${\cI}\subset \ker
\pi_E$, its dimension is less than $q$.

\subsection{Noncommutative pseudodifferential calculus}
Noncommutative pseudodifferential calculus for a smooth spectral
triple over an unital algebra $\cA$ was introduced by Connes and
Moscovici \cite{Co-M,Sp-view}. Their definition was extended to
the non-unital case in \cite{egorgeo}.

Assume that $({\mathcal A},{\mathcal H},D)$ is a $QC_0^\infty$
spectral triple. By the spectral theorem, for any $s\in\RR$, the
operator $\langle D\rangle^s=(D^2+I)^{s/2}$ is a well-defined
positive self-adjoint operator in $\cH$, which is unbounded for
$s>0$. For any $s\geq 0$, define by $\cH^s$ the domain of
$\langle D\rangle^s$, and, for $s<0$, put $\cH^s=(\cH^{-s})^*$.
Let also $\cH^{\infty}=\bigcap_{s\geq 0}\cH^s, \quad
\cH^{-\infty}=(\cH^{\infty})^*$.

\begin{defn}
We say that an operator $P$ in $\cH^{-\infty}$ belongs to the
class $\Psi^*_0({\mathcal A})$, if it admits an asymptotic
expansion:
\[ P\sim \sum_{j=0}^{+\infty}b_{q-j}\langle D \rangle^{q-j},
\quad b_{q-j}\in {\mathcal B},
\]
that means that, for any $N$,
\[
P - \left(b_q\langle D\rangle^q + b_{q-1}\langle
D\rangle^{q-1}+\ldots+b_{-N}\langle D\rangle^{-N}\right)\in {\rm
OP}_0^{-N-1}.
\]
\end{defn}

For the spectral triple $(\cA,\cH,\cD)$ associated with a compact
Riemannian manifold $(M,g)$, one can show that $\cH^s = H^s(M,E)$
for any $s$ and $\Psi^*_0({\mathcal A})=\Psi^0(M)$.

Let $(M,{\mathcal F})$ be a compact foliated manifold. Consider a
spectral triple $(\cA,\cH,D)$ described by (T1), (T2), (T3). One
can show that $H^s(M,E)\subset \cH^s$ for any $s\geq 0$ and $\cH^s
\subset H^s(M,E)$ for any $s<0$.

\begin{defn}\label{d:ideals}
The class $\cL^1(\cH^{-\infty},\cH^{\infty})$ consists of all
bounded operators $A$ in $\cH^{\infty}$ such that, for any real
$s$ and $r$, the operator $\langle D\rangle^rA\langle
D\rangle^{-s}$ extends to a trace class operator in $L^2(M,E)$.
\end{defn}

The class $\cL^1(\cH^{-\infty},\cH^{\infty})$ is an involutive
subalgebra in $\cL(\cH)$, and any operator with the smooth kernel
belongs to $\cL^1(\cH^{-\infty},\cH^{\infty})$.

\begin{prop}
(1) Any element $b\in \cB$ can be written as
\[
b=B+T, \quad B\in \Psi_{sc}^{0,-\infty}(M,\cF,E), \quad T\in
\cL^1(\cH^{-\infty},\cH^{\infty}).
\]
(2) The algebra $\Psi^*_0({\mathcal A})$ is contained in
$\Psi_{sc}^{*,-\infty}(M,{\mathcal F},E)+{\rm OP}_0^{-N}$ for any
$N$.
\end{prop}

\subsection{Noncommutative geodesic flow}
The definitions of the unitary co\-tan\-gent bundle and the
noncommutative geodesic flow associated with a $QC_0^\infty$
spectral triple $({\mathcal A},{\mathcal H},D)$ are motivated by
the relation (\ref{e:cosphere}) and the Egorov theorem,
Theorem~\ref{t:cl-egorov}.

Put ${\mathcal C}_0={\rm OP}_0^0\bigcap \Psi^*_0({\mathcal A})$.
Let $\bar{\mathcal C}_0$ be the closure of ${\mathcal C}_0$ in
${\mathcal L}(\cH)$. For any $T\in \cL(\cH)$, define
\begin{equation}\label{e:alphat}
\alpha_t(T)=e^{it\langle D\rangle}Te^{-it\langle D\rangle}, \quad
t\in {\RR}.
\end{equation}

\begin{defn}\cite{Sp-view,egorgeo}
The unitary cotangent bundle $S^*{\mathcal A}$ is defined as the
quotient of the $C^*$-algebra, generated by the union of all
spaces of the form $\alpha_t(\bar{\mathcal C}_0)$ with $t\in\RR$
and $\mathcal K$, by its ideal ${\mathcal K}$.
\end{defn}

\begin{defn}\cite{Sp-view,egorgeo}
The noncommutative geodesic flow is the one-parameter group
$\alpha_t$ of automorphisms of the algebra $S^*{\mathcal A}$
defined by (\ref{e:alphat}).
\end{defn}

As shown in \cite{Sp-view}, for the spectral triple
$(\cA,\cH,\cD)$ associated with a compact Riemannian manifold
$(M,g)$, the unitary cotangent bundle $S^*{\mathcal A}$ is the
algebra $C(S^*M)$ of continuous functions on the cosphere bundle
$S^*M$ and the noncommutative geodesic flow on $S^*{\mathcal A}$
is induced by the restriction of the geodesic flow to $S^*M$.

Theorem~\ref{Egorov} allows to give a description of the
noncommutative flow defined by a spectral triple associated with a
Riemannian foliation in the case when $E$ is the trivial line
bundle (see \cite{egorgeo}).

\begin{thm}
\label{noncom:flow} Consider a spectral triple $({\mathcal
A},{\mathcal H},D)$ defined in (T1), (T2), (T3) when $E$ is the
trivial line bundle and the subprincipal symbol of $D^2$ vanishes.
There is a nontrivial $\ast$-homomorphism $P : S^*{\mathcal
A}\rightarrow \bar{S}^{0}(G_{{\mathcal F}_N}, |T{\mathcal
G}_N|^{1/2})$ such that the following diagram commutes:
\[
  \begin{CD}
S^*{\mathcal A} @>\alpha_t>> S^*{\mathcal A}\\ @VPVV @VVPV
\\ \bar{S}^{0}(G_{{\mathcal F}_N}, |T{\mathcal
G}_N|^{1/2}) @>F^*_t>> \bar{S}^{0}(G_{{\mathcal F}_N},
|T{\mathcal G}_N|^{1/2})
  \end{CD}
\]
Here $F^*_t$ is the transverse Hamiltonian flow on
$C^{\infty}_{prop}(G_{{\mathcal F}_N},|T{\mathcal G}_N|^{1/2})$
associated with $\sqrt{a_2}$, where $a_2\in S^2(\tilde{T}^*M)$ is
the principal symbol of $D^2$.
\end{thm}

An extension of this theorem to the case of an arbitrary vector
bundle $E$ is directly related with an answer to
Problem~\ref{p:egorov}. The $\ast$-homomorphism $P$ is essentially
induced by the principal symbol map $\bar{\sigma}$. Therefore, a
more precise information on injectivity and surjectivity
properties of $P$ depends on answers to Questions~\ref{q:seq1} and
\ref{q:seq2}.

\end{document}